\documentclass[11pt]{amsart}
\usepackage{geometry}                
\usepackage{graphicx}
\usepackage{amssymb}
\usepackage{amsmath,mathrsfs}
\usepackage{amsfonts}
\usepackage{bbold,bm}
\usepackage{epstopdf}
\usepackage{color}
\usepackage{natbib}
\usepackage{algorithmicx,algorithm,algpseudocode}
\newtheorem{theorem}{Theorem}[section]

\newtheorem{remark}[theorem]{Remark}

\title{Discrete gradients for computational Bayesian inference}
\author{Sahani Pathiraja \and Sebastian Reich}

\begin{document}
\maketitle

\begin{abstract} 
In this paper, we exploit the gradient flow structure of continuous-time formulations of 
Bayesian inference in terms of their numerical time-stepping. We focus on two particular examples,
namely, the continuous-time ensemble Kalman--Bucy filter and a particle discretisation of the Fokker--Planck equation
associated to Brownian dynamics. Both formulations can lead to stiff differential equations which require special
numerical methods for their efficient numerical implementation. We compare discrete gradient methods to alternative
semi-implicit and other iterative implementations of the underlying Bayesian inference problems. 
\end{abstract}

%
%
\section{Introduction}
%
%

A number of algorithmic approaches to Bayesian inference \citep{sr:Robert2001,sr:Tarantola,sr:kaipio} can be rephrased in terms of
evolution equations of the form
\begin{equation} \label{eq:ODE1}
{\rm d}z_\tau = -A(z_\tau)\nabla_z V(z_\tau)\,{\rm d}\tau ,
\end{equation}
$\tau \ge 0$, with $A(z)$ being symmetric positive semi-definite and $V:\mathbb{R}^{N_z} \to \mathbb{R}$ being an appropriate
potential. It holds that
\begin{equation}\label{eq:decay}
{\rm d}V(z_\tau) = -\nabla_z V(z_\tau) \cdot A(z_\tau) \nabla_z V(z_\tau)\,{\rm d}\tau \le 0
\end{equation}
along solutions $z_\tau \in \mathbb{R}^{N_z}$ of (\ref{eq:ODE1}). ODEs of the form (\ref{eq:ODE1}) are referred to as gradient dynamics or of gradient flow structure \citep{sr:SH96}. Generally speaking, numerical approximations of (\ref{eq:ODE1}) will not preserve 
the decay property (\ref{eq:decay}) unless the time-step $\Delta \tau>0$ is chosen sufficiently small. Such time-step restrictions
can be problematic if the ODE system (\ref{eq:ODE1}) is stiff \citep{sr:ascher08}.  Discrete gradient methods \citep{sr:McLQB99} are attractive in this context, as they have been devised such that (\ref{eq:decay}) holds regardless of the chosen step-size $\Delta \tau$. 
However, these methods are implicit, which implies an increased computational cost per time-step. In this paper, we explore the 
pros and cons of discrete gradient methods in the context of computational Bayesian inference.
More specifically, we consider two particular examples of such gradient dynamics arising from particle methods 
for sequential data assimilation \citep{sr:stuart15,sr:reichcotter15} and from a transformation approach to posterior sampling 
\citep{sr:crisan10,sr:daum11,sr:reich10}. In both cases, the decay property (\ref{eq:decay}) plays a crucial role in
ensuring the statistical consistency of the associated inference method when the state variable $z_\tau$ is treated as a random variable 
with the initial $z_0$ following a given prior distribution $\pi_0$.

The paper is structured as follows. Section \ref{sec:examples} provides the necessary background for the two
considered continuous-time gradient flow systems, namely, the ensemble Kalman--Bucy filter (EnKBF) and
a particle discretisation of the Fokker--Planck equations associated to Brownian dynamics. Section \ref{sec:DG}
summarises the discrete gradient method and its numerical implementation. We also discuss the semi-implicit Euler method as
an alternative time-stepping method. Numerical results and a comparison between different methods is provided in 
Section \ref{sec:Numerics}. The paper closes with a summary in Section \ref{sec:Summary}.

%
%
\section{Examples of gradient systems for Bayesian inference} \label{sec:examples}
%
%

The following two examples of gradient flow systems are both concerned with dynamically transforming a set of Monte Carlo 
samples $\{x^i_0\}$, $i=1,\ldots,M$, from a prior distribution $\pi_0(x)$ into samples $\{x^i_\tau\}$ from a posterior distribution
\begin{equation}
\label{eq:post}
\pi_\ast(x) \propto \pi(y|x)\,\pi_0(x),
\end{equation}
which can be used to approximate expectation values of a function $g$ 
of a random variable $X_\ast$ with distribution given by (\ref{eq:post}), that is:
\begin{equation} \label{eq:EV}
\mathbb{E} [g(X_\ast)] \approx \frac{1}{M} \sum_{i=1}^M g(x_\tau^i) .
\end{equation}
The first example, namely, the ensemble Kalman--Bucy filter \citep{sr:br10,sr:br10b,sr:br11}, achieves this for $\tau=1$, while the second example, based
on a particle discretisation of the Fokker--Planck equation \citep{sr:R90,sr:DM90,sr:reich19}, requires $\tau\to \infty$. We now describe both approaches in
more detail.

\subsection{Ensemble Kalman--Bucy filter}

We have $M$ particles $x_\tau^i \in \mathbb{R}^{N_x}$ that evolve in time according to the gradient dynamics \citep{sr:br10,sr:br10b}
\begin{equation} \label{eq:EnKBF}
{\rm d} x_\tau^i = -P_\tau^{xx} \nabla_{x^i} V(\{x_\tau^j\})\,{\rm d}\tau, \qquad i=1,\ldots,M,
\end{equation}
with the symmetric positive semi-definite (covariance) matrix $P_\tau^{xx}$ given by
\begin{equation}
P_\tau^{xx} = \frac{1}{M-1} \sum_{j=1}^M (x_\tau^j-\overline{x}_\tau)(x_\tau^j-\overline{x}_\tau)^{\rm T}, \qquad
\overline{x}_\tau = \frac{1}{M}\sum_{j=1}^M x_\tau^j
\end{equation}
and the potential $V:\mathbb{R}^{N_x\times M}\to \mathbb{R}$ by\footnote{Different variants of $V$ can be used when the EnKBF is
used for minimisation only. See, for example, \cite{sr:KS18,sr:reich19}.}
\begin{equation} \label{eq:potential1}
V(\{x_\tau^j\}) = \frac{M}{2}\left\{ S(\overline{x}_\tau) + \frac{1}{M}\sum_{i=1}^M S(x_\tau^i)\right\}.
\end{equation}
Here $S(x)$ denotes the negative log-likelihood function, that is,
\begin{equation}
S(x) = -\log \pi(y|x).
\end{equation}
We are primarily interested in data models of the form
\begin{equation} \label{eq:obs}
y = h(x) + R^{1/2} \Xi, \qquad \Xi \sim {\rm N}(0,I),
\end{equation}
and, hence,
\begin{equation} \label{eq:S}
S(x) = \frac{1}{2}(h(x)-y)^{\rm T}R^{-1}(h(x)-y).
\end{equation}
The ODE system (\ref{eq:EnKBF}) is typically solved over a unit time interval, that is, $\tau\in [0,1]$. The initial particles $\{x_0^i\}$ 
represents the prior distribution while the particles $\{x_\tau^i\}$ at final time approximate the posterior distribution $\pi_\ast$ and can be used
in (\ref{eq:EV}) with $\tau=1$.

The ODE system (\ref{eq:EnKBF}) is of the form (\ref{eq:ODE1}) with 
\begin{equation} \label{eq:Z}
z_\tau = \begin{pmatrix} x_\tau^1 \\ x_\tau^2 \\ \vdots \\ x_\tau^M\end{pmatrix}
\end{equation}
and $A(z_\tau)$ being a block diagonal matrix with repeated entries $P_\tau^{xx}$. Numerical problems arising from the stiffness of
the gradient dynamics have been discussed by \cite{sr:akir11}. 

We mention that there are derivative-free extensions of the EnKBF to non-linear forward maps $h$ \citep{sr:br11,sr:reichcotter15}. 
Consider, for example,
\begin{equation}
{\rm d}x_\tau^i = -P_\tau^{xh}R^{-1}\left(\frac{1}{2}(h(x_\tau^i)+\overline{h}_\tau) - y\right)
\end{equation}
with
\begin{equation}
P_\tau^{xh} = \frac{1}{M-1} \sum_{j=1}^M (x_\tau^j-\overline{x}_\tau)(h(x_\tau^j)-\overline{h}_\tau)^{\rm T}, \qquad
\overline{h}_\tau = \frac{1}{M}\sum_{j=1}^M h(x_\tau^j).
\end{equation}
This formulation is no longer of gradient flow structure unless $h$ is linear. There
are also stochastic formulations of the EnKBF \citep{sr:br11,sr:stuart15,sr:reichcotter15}.

We finally mention that EnKBF is closely related to iterative implementations of the ensemble Kalman filter as
considered, for example, by \cite{sr:reynolds12,sr:sakov12,sr:oliver13}. Those implementations can also be thought of
particular time-stepping methods for the continuous-time EnKBF. We will return to this aspect in Section \ref{sec:EnKBF2}.
Iterative implementations of the Kalman--Bucy filter are also related to the natural gradient approach from statistical learning 
\citep{sr:O17}.

\subsection{Particle flow Fokker--Planck dynamics} \label{sec:FPD}

Given $M$ particles $x_\tau^i$ and a kernel function $\psi(x)\ge 0$ of a Reproducing Kernel Hilbert Space (RKHS) $\mathcal{H}$, 
we approximate their density in $\mathcal{H}$ by
\begin{equation} \label{eq:PDF_H}
\widetilde{\pi}_\tau(x) = \frac{1}{M} \sum_{j=1}^M \psi(x-x_\tau^j) = \int \psi(x-x') \pi_\tau(x'){\rm d}x',
\end{equation}
where $\pi_\tau$ stands for the empirical measure
\begin{equation}
\pi_\tau (x) = \frac{1}{M}\sum_{i=1}^M \delta (x-x_\tau^i)
\end{equation}
and we have assumed that $\psi(x)=\psi(-x)$ as well as
\begin{equation}
\int \psi(x)\,{\rm d}x = 1.
\end{equation}

Let us denote the standard $L_2$-inner product of two functions $f$ and $g$ by $\langle f,g\rangle$ and the associated inner product
in $\mathcal{H}$ by $\langle f,g\rangle_\mathcal{H}$. The Kullback--Leibler divergence between $\widetilde{\pi}_\tau(x)$ and the posterior 
target PDF $\pi_\ast(x):=\pi(x|y)$ is then 
approximated in $\mathcal{H}$ by
\begin{align}
V(\{x_\tau^j\}) &= \Bigl\langle \widetilde{\pi}_\tau, \log \left(\frac{\widetilde{\pi}_\tau}{\pi_\ast} \right) \Bigr\rangle_{\mathcal{H}} = \Bigl\langle 
\pi_\tau , \log \left( \frac{\widetilde{\pi}_\tau}{\pi_\ast} \right) \Bigr\rangle \\
& = \frac{1}{M} \sum_{j=1}^M \left\{\log \widetilde{\pi}_\tau(x_\tau^j) - \log \pi_\ast(x_\tau^j)\right\} \\ & = \frac{1}{M} \sum_{j=1}^M \left\{
\log \left( \frac{1}{M} \sum_{l=1}^M \psi(x_\tau^j-x_\tau^l) \right) - \log \pi_\ast(x_\tau^j) \right\}
\end{align}
and the particle flow Fokker--Planck dynamics is given by 
\begin{equation} \label{eq:DFP1}
{\rm d}x_\tau^i = - M \nabla_{x^i} V(\{x_\tau^j\})\,{\rm d}\tau,
\end{equation}
for $i=1,\ldots,M$. See \cite{sr:reich19} for more details. The gradient of $V$ with respect to $x_\tau^i$ is given by
\begin{equation}
\nabla_{x_\tau^i} V(\{x_\tau^j\}) = \frac{1}{M} \left\{ \nabla_{x_\tau^i} \log \widetilde{\pi}_\tau (x_\tau^i) + 
\frac{1}{M} \sum_{j\not= i} \frac{1}{\widetilde{\pi}_\tau (x_\tau^j)}\nabla_{x_\tau^i} \psi(x_\tau^i-x_\tau^j) -
\nabla_{x_\tau^i} \log \pi_\ast(x_\tau^i)  \right\}.
\end{equation}
We note that one can replace the normalised kernel $\psi$ by an unnormalised version $\widetilde{\psi}$ without changing the
gradient dynamics. Furthermore, the ODE system (\ref{eq:DFP1}) is of the form (\ref{eq:ODE1}) with $z_\tau$ defined by 
(\ref{eq:Z}) and $A(z_\tau) = M I$ where $I$ is the identity matrix.

Let us look at the evolution equations (\ref{eq:DFP1}) and their geometric properties in some more detail.
The variational derivative of the Kullback--Leibler divergence in $\mathcal{H}$, 
\begin{equation}
{\rm KL}_\mathcal{H}(\pi_\tau|\pi^\ast) = \Bigl\langle \pi_\tau, \log \left( \frac{\widetilde{\pi}_\tau}{\pi_\ast} \right) \Bigr\rangle,
\end{equation}
with respect to $\pi_\tau$ is given by
\begin{equation} \label{eq:VD}
\frac{\delta {\rm KL}_\mathcal{H}(\pi_\tau|\pi^\ast)}{\delta \pi_\tau}  (x) =
\log \widetilde{\pi}_\tau (x) - \log \pi_\ast (x) + \int \psi(x-x') \frac{\pi_\tau (x')}{\widetilde{\pi}_\tau (x')} {\rm d}x'
\end{equation}
and, hence, one finds that
\begin{equation}
M \nabla_{x_\tau^i} V(\{x_\tau^j\}) = \nabla_x \left(\frac{\delta {\rm KL}_\mathcal{H}(\pi_\tau|\pi^\ast)}{\delta \pi_\tau} \right) (x_\tau^i) .
\end{equation}
Furthermore, as $\tau \to \infty$, the particle flow Fokker--Planck dynamics (\ref{eq:DFP1}) leads to
\begin{equation}
\nabla_x \left(  \frac{\delta {\rm KL}_\mathcal{H}(\pi_\infty |\pi^\ast)}{\delta \pi_\infty} \right) (x_\infty^i) = 0
\end{equation}
and the asymptotic particle positions $\{x_\infty^i\}$ can be used to approximate expectation values 
with respect to $X_\ast \sim \pi_\ast$ based on
\begin{equation}
\pi^\ast (x) \propto \widetilde{\pi}_\infty (x) \,e^{l(x)}, \qquad l(x) =  \int \psi(x-x') \frac{\pi_\infty (x')}{\widetilde{\pi}_\infty (x')} {\rm d}x' ,
\end{equation}
which follows from (\ref{eq:VD}) for an equilibrium density $\widetilde{\pi}_\infty$, that is, 
\begin{equation}
\frac{\delta {\rm KL}_\mathcal{H}(\pi_\infty|\pi^\ast)}{\delta \pi_\tau}  (x) = \mbox{const.}
\end{equation}
We use the simpler, less accurate in the limit $M\to \infty$ empirical 
approximation (\ref{eq:EV}) for $\tau>0$ sufficiently large in the numerical experiments of Section \ref{sec:Numerics}.

For $M\to \infty$ and the kernel function $\psi (x)$ approaching a Dirac delta function,
the evolution equations (\ref{eq:DFP1}) provide a mesh-free particle approximation to the Fokker--Planck 
equation for the marginal PDFs $p_t$ of Brownian dynamics
\begin{equation} \label{eq:BD1}
{\rm d}X_\tau = \nabla_x \log \pi_\ast(X_\tau)\,{\rm d}\tau + \sqrt{2} {\rm d}W_\tau
\end{equation}
\citep{sr:reichcotter15}. It holds under suitable conditions on $\pi_\ast$ that the marginal distributions $p_\tau$ of the solutions $X_\tau$ to the SDE (\ref{eq:BD1}) 
satisfy
\begin{equation}
\lim_{\tau \to \infty} p_\tau = \pi_\ast 
\end{equation}
in a weak sense \citep{sr:P14}.

As a specific example, let us consider the RKHS $\mathcal{H}$ with Gaussian kernel
\begin{equation} \label{eq:kernel}
\psi(x-x') := {\rm n}(x;x',B) \propto \widetilde{\psi}(x-x') := \exp \left( -\frac{1}{2} (x-x')^{\rm T}B^{-1}(x-x') \right)
\end{equation}
for given covariance matrix $B$. Here ${\rm n}(x,\overline{x},B)$ denotes the probability density function of a Gaussian
random variable with mean $\overline{x}$ and covariance matrix $B$. We now provide a construction for the required
covariance matrix $B$ and the initial particle positions $x_0^i$ that yields a kernel approximation to the given 
prior distribution $\pi_0$. Given a set of $M$ samples $\widehat{x}^i_0$ from the prior distribution 
$\pi_0$, the associated empirical measure is given by
\begin{equation}
\pi_{\rm em}(x) = \frac{1}{M} \sum_{i=1}^M \delta (x-\widehat{x}_0^i)
\end{equation}
with empirical mean $\overline{x}_0$ and empirical covariance matrix $P^{xx}_0$. Setting 
\begin{equation} \label{eq:shrink}
B = (2\alpha - \alpha^2) P_0^{xx}
\end{equation}
for $\alpha \in (0,1]$ in (\ref{eq:kernel}), 
the associated measure (\ref{eq:PDF_H}) in $\mathcal{H}$ is defined as follows:
\begin{equation} \label{eq:GMFP}
\widetilde{\pi}_0(x) = \frac{1}{M} \sum_{i=1}^M {\rm n}(x;x^i_0,B)
\end{equation}
with particle positions 
\begin{equation}
x^i_0 = \widehat{x}_0^i - \alpha(\widehat{x}_0^i - \overline{x}_0).
\end{equation} 
In other words, the covariance matrix $B$ and the particle positions $x^i_0$ are chosen such that 
the mean and covariance matrix under $\widetilde{\pi}_0$ are 
identical to $\overline{x}_0$ and $P_0^{xx}$, respectively, that is,
\begin{equation}
\sum_{i=1}^M x_0^i = \sum_{i=1}^M \widehat{x}^i_0
\end{equation}
and
\begin{equation}
P_0^{xx} = (2\alpha - \alpha^2) P_0^{xx} + \frac{1}{M-1}\sum_{i=1}^M (x^i_0-\overline{x}_0)(x^i_0-\overline{x}_0)^{\rm T}.
\end{equation}
The initial particle positions $\{x_0^i\}$ are now evolved under the ODEs (\ref{eq:DFP1}). Note that $\alpha=1$ leads to $x^i_0 = \overline{x}_0$ 
and points of $\nabla_{x^i} V = 0$, that is, stationary points of (\ref{eq:DFP1}) agree with critical points of the posterior distribution $\pi_\ast$. 
It is possible to extend the dynamics in the particles $\{x_\tau^i\}$ by an evolution equation in the kernel matrix $B$ \citep{sr:Y97}. 

Furthermore, one can put (\ref{eq:DFP1}) into the more general form of (\ref{eq:ODE1}), that is,
\begin{equation}\label{eq:DFP1b}
{\rm d}x_\tau^i = -\mathcal{A}(\{x_\tau^j\}) \nabla_{x^i} V(\{x_\tau^j\})\,{\rm d}\tau,
\end{equation}
for $i=1,\ldots,M$. A particular choice is given by $\mathcal{A}(\{x_\tau^j\}) = M\,P_t^{xx}$. A related preconditioned Brownian dynamics formulation 
(\ref{eq:BD1}) has been proposed recently by \cite{sr:GHLS19}.

We finally mention that the Stein variational approach \citep{sr:LW16} has recently become popular as an alternative dynamical approach for transforming a set of Monte Carlo samples from a prior distribution into samples from a posterior distribution. However, it is currently 
unknown whether the resulting interacting particle systems possess a gradient flow structure. See \cite{sr:DCSMS18} for a preconditioned 
Stein variational formulation in the spirit of (\ref{eq:DFP1b}).

%
%
\section{Discrete gradient time-stepping} \label{sec:DG}
%
%

Discrete gradient methods are numerical time-stepping methods for ODEs of type (\ref{eq:ODE1}) that allow one to maintain the
decay property (\ref{eq:decay}) under arbitrary choices of the step-size $\Delta \tau>0$. See  \cite{sr:McLQB99} for an overview of
such methods. Here we focus on the following particular form of discrete gradient dynamics:
\begin{equation} \label{eq:DG1}
z_{n+1}-z_n = -\Delta \tau A(z_{n+\theta}) \overline{\nabla}_z V(z_{n+1},z_n), \qquad z_{n+\theta} = \theta z_{n+1}+ (1-\theta) z_n,
\end{equation}
$\theta \in [0,1]$, with the discrete gradient
\begin{equation}
\overline{\nabla}_z V(z_{n+1},z_n) := \frac{V(z_{n+1})-V(z_n)}{\nabla_z V(z_{n+\theta}) \cdot (z_{n+1}-z_n)} \nabla_z V(z_{n+\theta})
\end{equation}
and symmetric positive definite matrix $A(z)$. It follows that
\begin{equation}
V(z_{n+1})-V(z_n) = 
- \Delta \tau \overline{\nabla}_z V(z_{n+1},z_n) \cdot A(z_{n+\theta}) \overline{\nabla}_z V(z_{n+1},z_n) \le 0, \label{eq:decay2}
\end{equation}
as desired. Note that (\ref{eq:decay2}) is a rather natural discrete representation of the continuous (\ref{eq:decay}).

There exists several variants of the method, for example, the time-symmetric formulation of \cite{sr:Gonzalez96} (that is, with $\theta = 1/2$) and the local update version of \cite{sr:reich96}. 

We rewrite (\ref{eq:DG1}) in the form
\begin{equation} \label{eq:DG2}
z_{n+1} - z_n = -\gamma_n\Delta \tau A(z_{n+\theta}) \nabla_z V(z_{n+\theta}),
\end{equation}
with the scalar factor $\gamma_n$ defined by
\begin{equation}
\gamma_n = \frac{V(z_{n+1})-V(z_n)}{\nabla_z V(z_{n+\theta}) \cdot (z_{n+1}-z_n)} .
\end{equation}
Hence one can think of the discrete gradient method (\ref{eq:DG1}) as a $\theta$-method \citep{sr:SH96} with an implicitly determined
effective step-size
\begin{equation}
\Delta \tau_{n} = \gamma_{n}\Delta \tau.
\end{equation}

Formulation (\ref{eq:DG2}) with $\theta >0$ leads naturally to the following iterative solution procedure at each time-step $\tau_n$:
\begin{itemize}
\item Set $l=0$, $z_{n+\theta}^l = z_n$, and $\gamma_n^l = 1$.
\item For $l\ge 0$ solve till convergence the minimisation problem
\begin{equation} \label{eq:MP1}
z_{n+\theta}^{l+1} = \arg \min_z \left\{ \frac{(z-z_n)^{\rm T} A(z_{n+\theta}^l)^{-1} (z-z_n)}{2\theta}
+ \gamma_n^l\Delta \tau V(z)\right\}
\end{equation}
and set
\begin{align}
z_{n+1}^{l+1} &= \theta^{-1}(z_{n+\theta}^{l+1}-(1-\theta)z_n),\\
\gamma_n^{l+1} &=  \frac{V(z_{n+1}^{l+1})-V(z_n)}{\nabla_z V(z_{n+\theta}^{l+1}) \cdot (z_{n+1}^{l+1}-z_n)}.
\end{align}
\end{itemize}

The minimisation problem (\ref{eq:MP1}) can be solved by the Gauss--Newton method \citep{sr:wright99} when applied to the EnKBF
with data model (\ref{eq:obs}). Hence only first order derivatives of the nonlinear forward operator $h$ are required. Alternatively, one
can use a quasi-Newton method \citep{sr:wright99} directly on the nonlinear equations arising from the chosen time-step method.
In case $A(z)$ is only positive semi-definite, the inverse in (\ref{eq:MP1}) should be replaced by the pseudo-inverse and the update
is performed in the subspace spanned by the image of $A(z_{n+\theta}^l)$.

In the subsequent numerical examples, we will compare the discrete gradient formulation (\ref{eq:DG1}) to the standard semi-implicit
Euler method
\begin{equation} \label{eq:SI1}
z_{n+1} -z_n = -\Delta \tau A(z_n) \nabla V(z_{n+1}),
\end{equation}
which can also be reformulated as a minimisation problem, that is,
\begin{equation} \label{eq:MP-SI1}
z_{n+1} = \arg \min_z \left\{ \frac{(z-z_n)^{\rm T} A(z_{n})^{-1} (z-z_n)}{2}
+ \Delta \tau V(z)\right\}.
\end{equation}
We note that (\ref{eq:MP-SI1}) is structurally similar to the minimisation problem (\ref{eq:MP1}). However, the semi-implicit Euler method
(\ref{eq:SI1}) requires the solution of such a minimisation problem only once per time-step while the discrete gradient method leads to a
more elaborate implementation. Furthermore, iteration (\ref{eq:SI1}) does not obey a discrete decay property 
of the form (\ref{eq:decay2}).

While the discrete decay property (\ref{eq:decay2}) is sufficient to reach a critical point of (\ref{eq:DFP1}) at $\tau \to \infty$, a sufficiently accurate
approximation of the ensemble Kalman--Bucy filter equations is required in order to reproduce the posterior distribution. This statement is confirmed
by the numerical experiments in the following section. We note in this context that the discrete gradient dynamics (\ref{eq:DG1}) is first order
for all $\theta \in (0,1]$ except for $\theta = 1/2$ when the method becomes second order.

%
%
\section{Numerical results} \label{sec:Numerics}
%
%

We study the behaviour of the proposed gradient flow systems and their discretisation for a linear and a nonlinear
one-dimensional inference problem, as well as for a partially observed nonlinear three-dimensional system. We start with the ensemble Kalman--Bucy filter formulation.

\subsection{Ensemble Kalman--Bucy filter} \label{sec:EnKBF2}

\subsubsection{Linear problem}
We consider the scalar inference problem with Gaussian prior $X \sim {\rm N}(m_0,\sigma_0)$ and forward model
\begin{equation}
Y = X + \Xi, \qquad \Xi \sim {\rm N}(0,r).
\end{equation}
The posterior distribution is also Gaussian with mean
\begin{equation}
m_\ast = m_0 + K(y-m_0),
\end{equation}
Kalman gain matrix $K = \sigma_0/(\sigma_0 + r)$, and variance
\begin{equation}
\sigma_\ast = \sigma _0- K\sigma_0 .
\end{equation}
The experimental parameters are set to $m_0= 1/2$, $\sigma_0 = 1$, $r = 0.02$, $y= 0.1$. 

\begin{figure}
\begin{center}
\includegraphics[width=0.45\textwidth,trim = 0 0 0 0,clip]{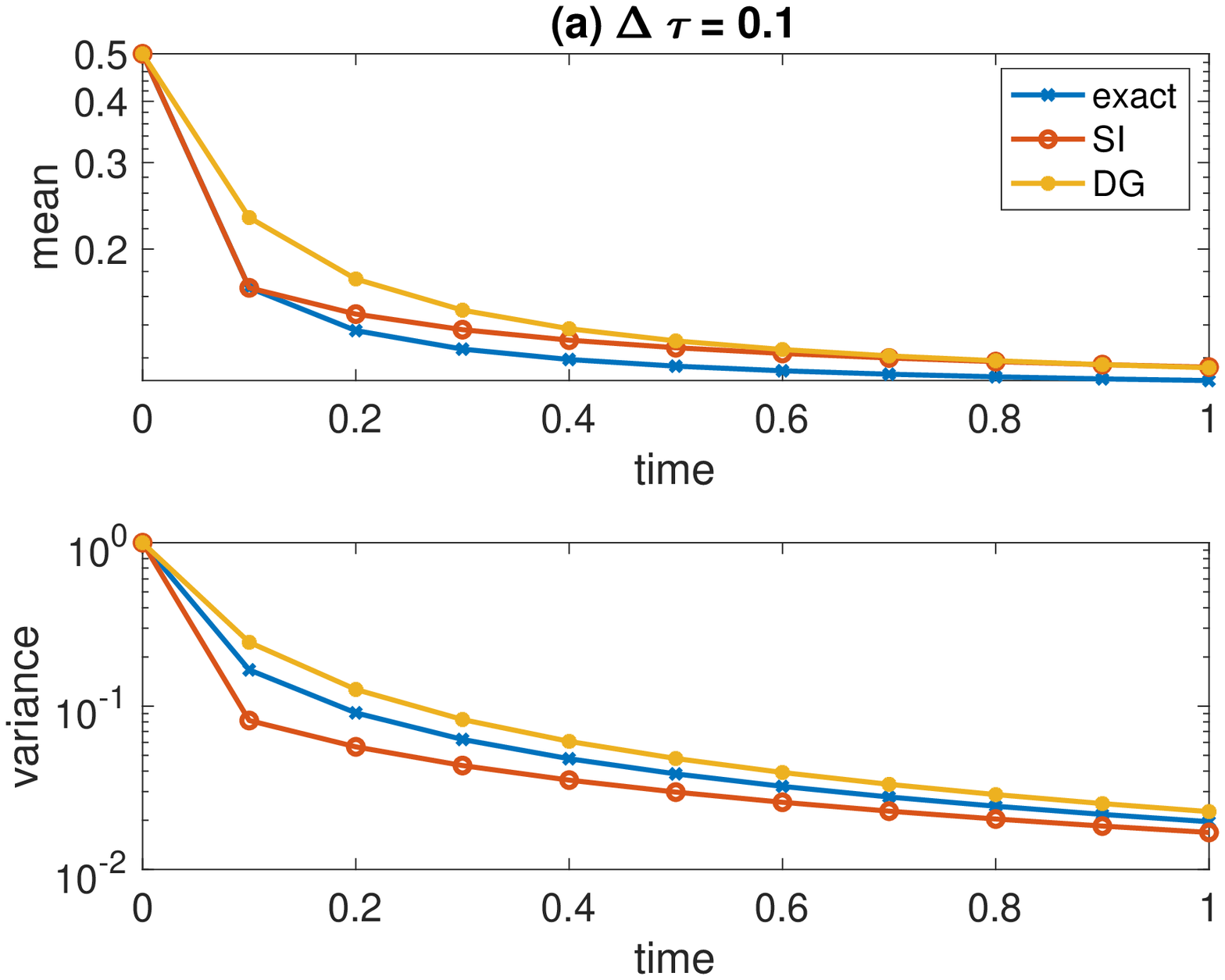} $\qquad$
\includegraphics[width=0.45\textwidth,trim = 0 0 0 0,clip]{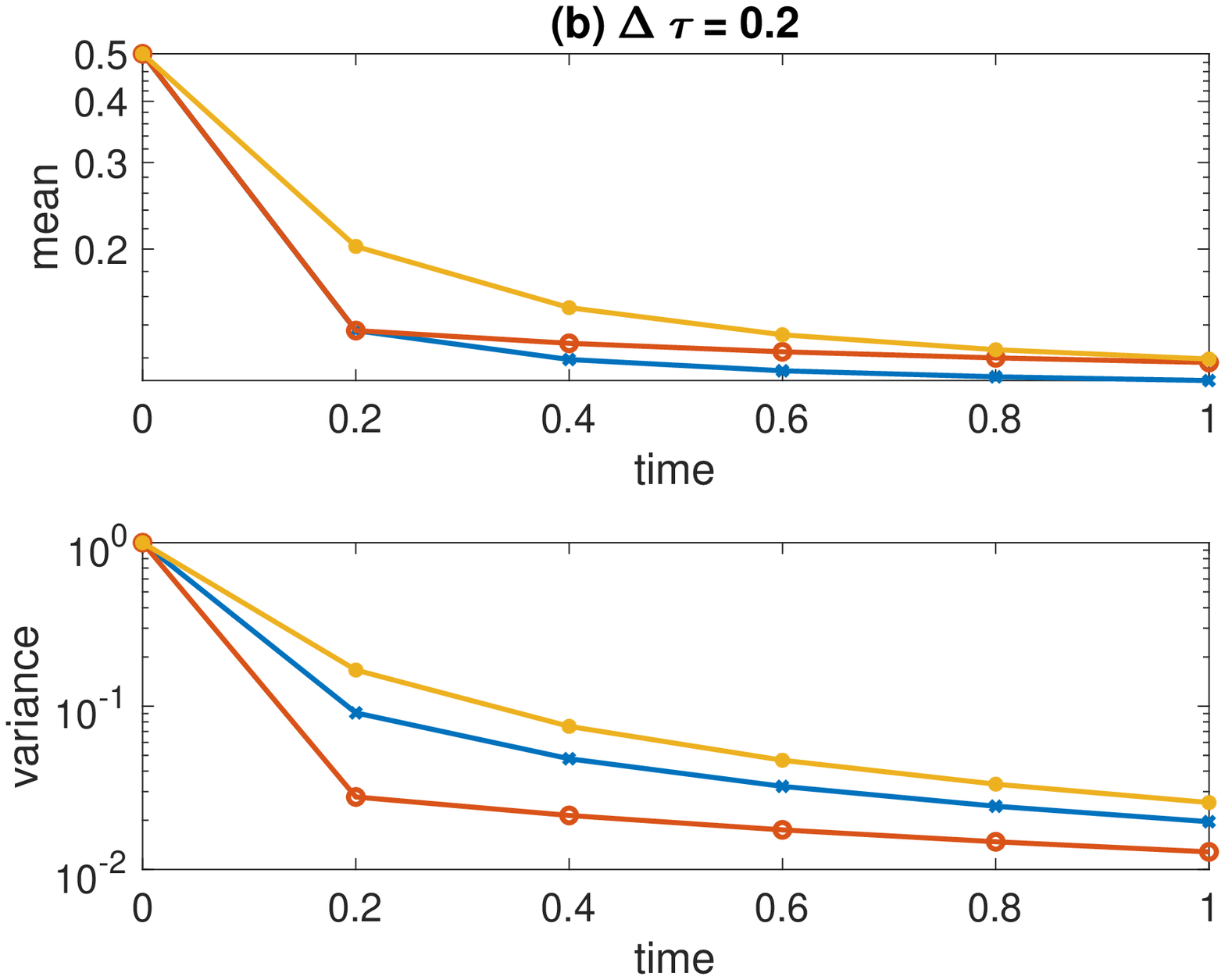} \\ \medskip
\includegraphics[width=0.45\textwidth,trim = 0 0 0 0,clip]{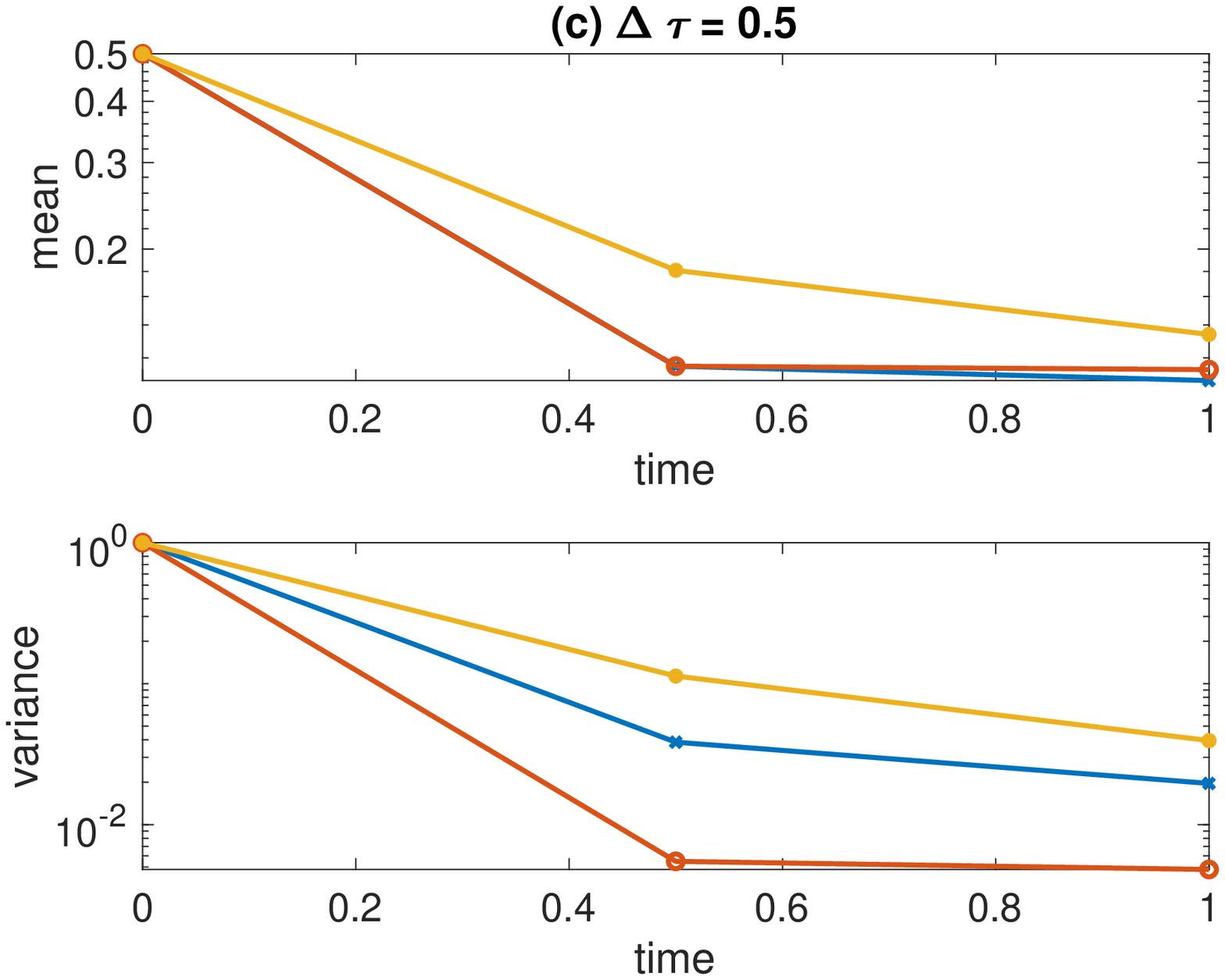} $\qquad$
\includegraphics[width=0.45\textwidth,trim = 0 0 0 0,clip]{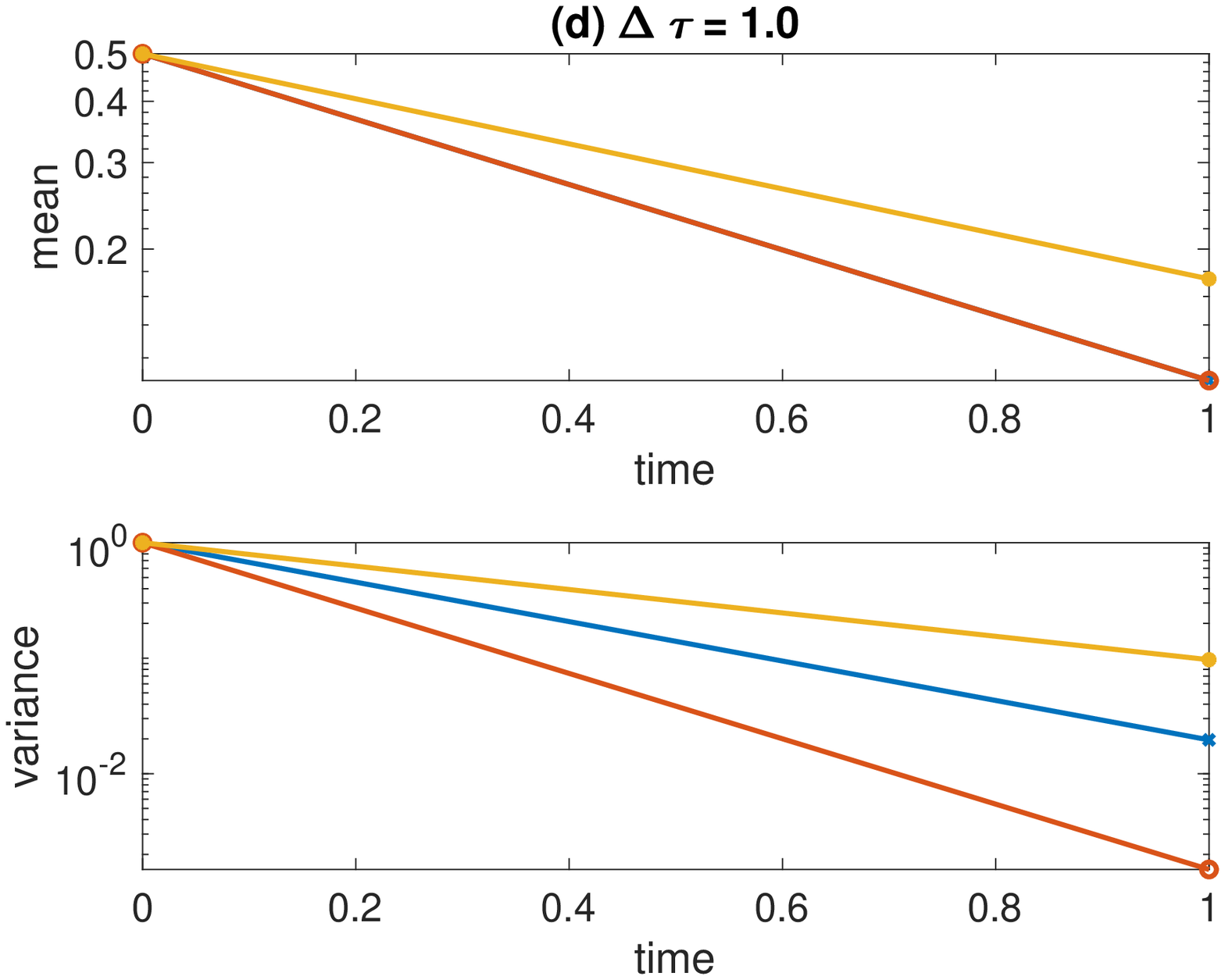} 
\end{center}
\caption{Exact solution and numerical approximations to the mean $m_\tau$ and variance $\sigma_\tau$ 
estimated by the EnKBF using the discrete gradient (DG) method with 
$\theta = 1$ and the semi-implicit (SI) Euler method with step-sizes (a) $\Delta \tau = 0.1$, (b) $\Delta \tau = 0.2$, (c) $\Delta \tau = 0.5$, (d) $\Delta \tau = 1.0$.}
\label{fig:figure1}
\end{figure}

We implement the EnKBF formulation with $M=2$ ensemble members. Their initial positions $x_0^i$ are uniquely determined by the
requirement that
\begin{equation}
m_0 = \frac{1}{2} (x_0^1 + x_0^2), \qquad \sigma_0 = (x_0^1-m)^2 + (x_0^2-m)^2 = \frac{1}{2} (x_0^2-x_0^1)^2.
\end{equation}
Furthermore, the analytic solutions $x_\tau^i$, $\tau\in (0,1]$, can be recovered from the associated solutions for the mean and 
variance, that is,
\begin{equation} \label{eq:exact_sol}
m_\tau = m - K_\tau(y-m), \qquad \sigma_\tau = \sigma - K_\tau \sigma,
\end{equation}
with time-dependent Kalman gain 
\begin{equation}
K_\tau = \frac{\sigma \tau}{\sigma \tau+ r}.
\end{equation}
The potential (\ref{eq:potential1}) reduces to
\begin{equation}
V(x^1,x^2) = \frac{1}{8r} \left\{ (x^1+x^2-2y)^2 + 2(x^1-y)^2 + 2(x^2-y)^2\right\}
\end{equation}
with gradients
\begin{equation}
\nabla_{x^1} V(x^1,x^2) = \frac{1}{4r} \left( 3x^1+x^2 \right) - \frac{y}{r} , \quad
\nabla_{x^2} V(x^1,x^2) = \frac{1}{4r} \left( x^1+3x^2 \right) - \frac{y}{r} .
\end{equation}
The matrix $A(x^1,x^2)$ is given by
\begin{equation}
A(x^1,x^2) = \frac{1}{2} \begin{pmatrix}  (x^2-x^1)^2 & 0 \\ 0 & (x^2-x^1)^2 \end{pmatrix}.
\end{equation}
We implement the discrete gradient method (\ref{eq:DG1}) with $\theta = 1$ and the semi-implicit Euler (\ref{eq:SI1}) method  
for $\Delta \tau \in \{0.1, 0.2, 0.5,1\}$. The relative small value of $r$ renders (\ref{eq:EnKBF}) stiff and the explicit Euler method is unstable for any
of the four chosen step-sizes $\Delta \tau$. The numerical results are compared to the analytic solutions as given by (\ref{eq:exact_sol}) and
can be found in Figure \ref{fig:figure1}. We find that the discrete gradient method systematically overestimates the variance while the semi-implicit
Euler method has the opposite effect on the variance. The semi-implicit Euler method also performs remarkably well with respect to the mean while
the discrete gradient method leads to relatively large errors in the final mean for step-sizes $\Delta \tau =0.5$ and $\Delta \tau = 1.0$.

\subsubsection{Nonlinear problem} \label{sec:numnonlinear}

We consider the following nonlinear forward map
\begin{equation}
h(x) = \frac{7}{12}x^3 - \frac{7}{2} x^2 + 8x
\end{equation}
and a measurement error $\Xi \sim {\rm N}(0,r)$ with $r=1$. The observed value is $y=2$ and the prior distribution is Gaussian with
$X\sim {\rm N}(-2,1/2)$. The posterior target distribution is given by
\begin{equation}
\pi^\ast(x) \propto e^{-(x+2)^2 - \frac{1}{2} (h(x)-2)^2}.
\end{equation}
The posterior mean and variance are given by $\overline{x}^\ast \approx 0.2095$ and
$\sigma^\ast \approx 0.0211$, respectively. See Example 5.9 in \cite{sr:reichcotter15} for more details. 

Given $M$ particles $x_n^i$ and a state vector
$z_n = (x_n^1,\ldots,x_n^M)^{\rm T}$ the required solution $z_{n+\theta}^{l+1} \in \mathbb{R}^{M}$ to the minimisation problem 
(\ref{eq:MP1}) can be recast in the  following form:
\begin{equation} \label{eq:MP2}
z_{n+\theta}^{l+1} = \arg \min_z \|f(z)\|_{C^l}^2,
\end{equation}
where $f:\mathbb{R}^{M} \to \mathbb{R}^{2M+1}$ is defined by
\begin{equation}
f(z) = \begin{pmatrix} z-z_n\\ h(\overline{x})-y \\ h(x^1)-y\\ \vdots \\ h(x^M)-y \end{pmatrix}, \quad 
z = \begin{pmatrix} x^1\\ x^2 \\ \vdots \\ x^M \end{pmatrix}, \quad \overline{x} = \frac{1}{M} \sum_{i=1}^M x^i,
\end{equation}
and the norm $\|\cdot\|_{C^l}$ by
\begin{equation}
\|f(z)\|_{C^l} = \frac{1}{2} f(z)^{\rm T} C^l f(z), \quad 
C^l = \begin{pmatrix} \frac{1}{\sigma_{n+\theta}^l} I & 0 & 0\\ 0 & \frac{\gamma_{n}^l \theta \Delta t M}{r} & 0 \\ 0 & 0 & 
\frac{\gamma_{n}^l \theta \Delta t}{r} I \end{pmatrix} .
\end{equation}
The Gauss--Newton method replaces the nonlinear minimisation problem
(\ref{eq:MP2}) by the linearised problem
\begin{equation} \label{eq:MP3}
z_{n+\theta}^{l+1} = \arg \min_z \|f(z_{n+\theta}^l) + Df(z_{n+\theta}^l) (z-z_{n+\theta}^l)\|_{C^l}^2.
\end{equation}
Here the Jacobian matrix $Df(z) \in \mathbb{R}^{2M+1\times M}$ is provided by
\begin{equation}
Df(z) = \begin{pmatrix} I \\ M^{-1} h'(\overline{x}) \mathbb{1} \\ D(c) \end{pmatrix}
\end{equation}
where $h'(x) = 7x^2/4 - 7x + 8$,
\begin{equation}
\mathbb{1} = \begin{pmatrix} 1 & 1 &\cdots & 1 \end{pmatrix} \in \mathbb{R}^{1\times M}, \qquad  
c = \begin{pmatrix} h'(x^1) & h'(x^2) &\cdots & h(x^M) \end{pmatrix} \in \mathbb{R}^{1\times M},
\end{equation}
and $D(c) \in \mathbb{R}^{M\times M}$ is a diagonal matrix with the entries of $c$ on its diagonal.

The semi-implicit Euler method (\ref{eq:SI1}) leads to a rather similar Gauss--Newton iteration
\begin{equation} \label{eq:MP-SI2}
z_{n+1}^{l+1} = \arg \min_z \|f(z_{n+1}^l) + Df(z_{n+1}^l) (z-z_{n+1}^l)\|_{C}^2
\end{equation}
with $C^l$ in (\ref{eq:MP2}) being replaced by
\begin{equation}
C = \begin{pmatrix} \frac{1}{\sigma_n} I & 0 & 0\\ 0 & \frac{\Delta t M}{r} & 0 \\ 0 & 0 & 
\frac{\Delta t}{r} I \end{pmatrix} 
\end{equation}
and $\theta = 1$.

We compare the results from the semi-implicit Euler and discrete gradient method with $\theta = 1$ to the 
gradient-free, explicit time-stepping method (IEnKF)
\begin{equation} \label{eq:GFE}
x_{n+1}^i - x_n^i = -\Delta t P_n^{xh}(\Delta t P_n^{hh}+R)^{-1}\left(\frac{1}{2}(h(x_n^i)+\overline{h}_n) - y\right),
\end{equation}
which has been discussed, for example, by \cite{sr:akir11,sr:BSW18,sr:dWRS18}. The discretisation (\ref{eq:GFE})
is also related to the DEnKF formulation of the ensemble Kalman as proposed by \cite{sr:sakov08}. 

\begin{figure}
\begin{center}
\includegraphics[width=0.45\textwidth,trim = 0 0 0 0,clip]{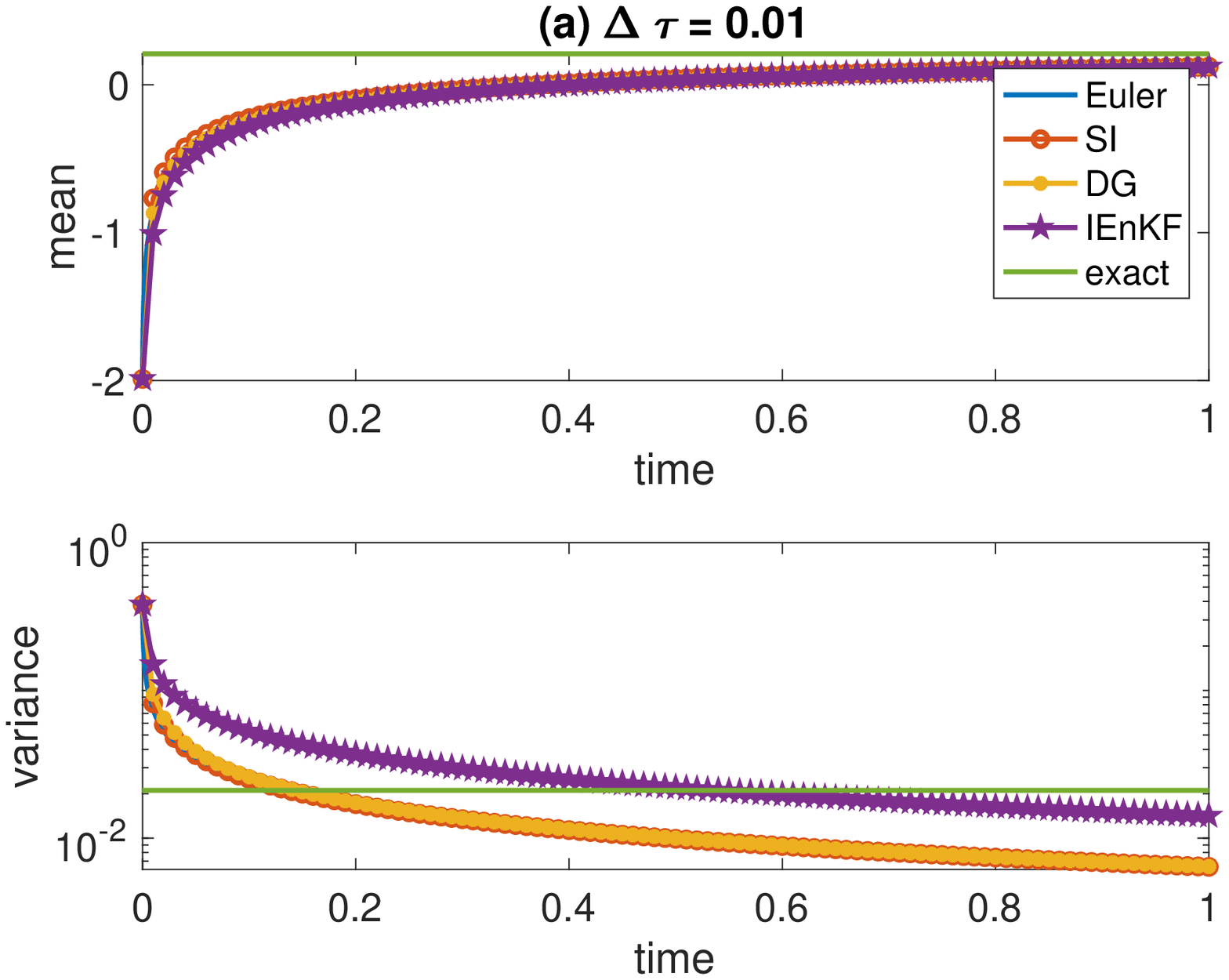} $\qquad$
\includegraphics[width=0.45\textwidth,trim = 0 0 0 0,clip]{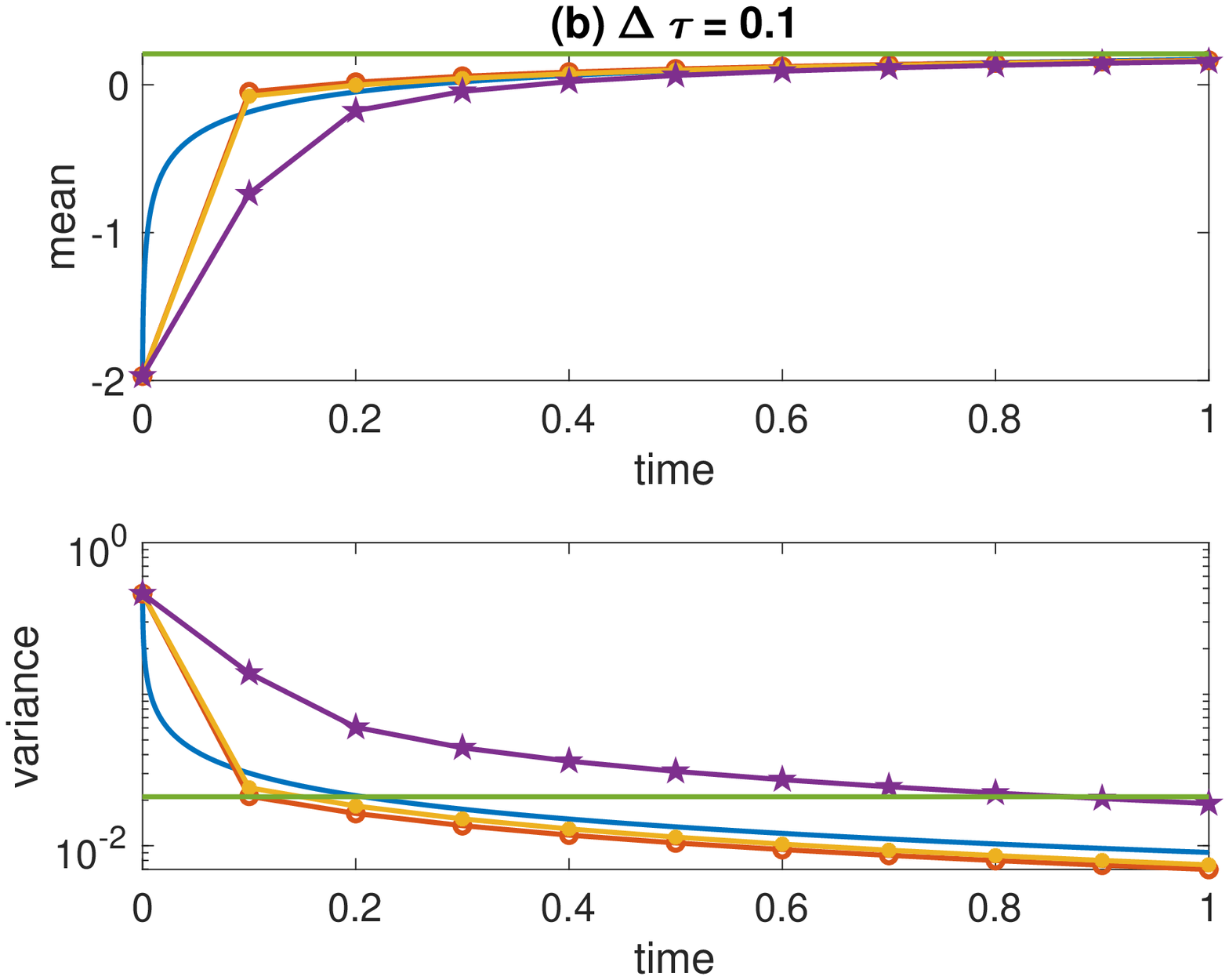} \\ \medskip
\includegraphics[width=0.45\textwidth,trim = 0 0 0 0,clip]{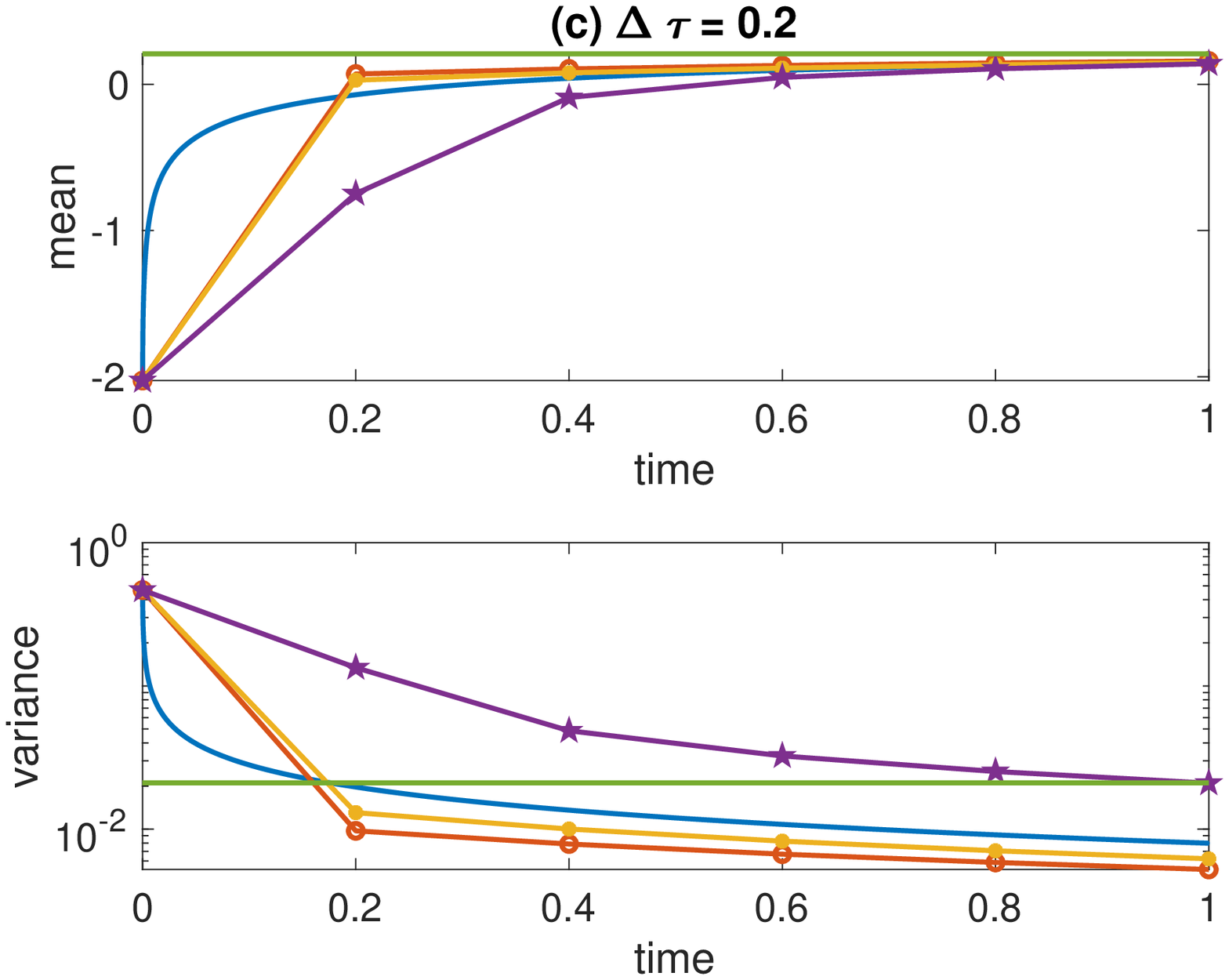} $\qquad$
\includegraphics[width=0.45\textwidth,trim = 0 0 0 0,clip]{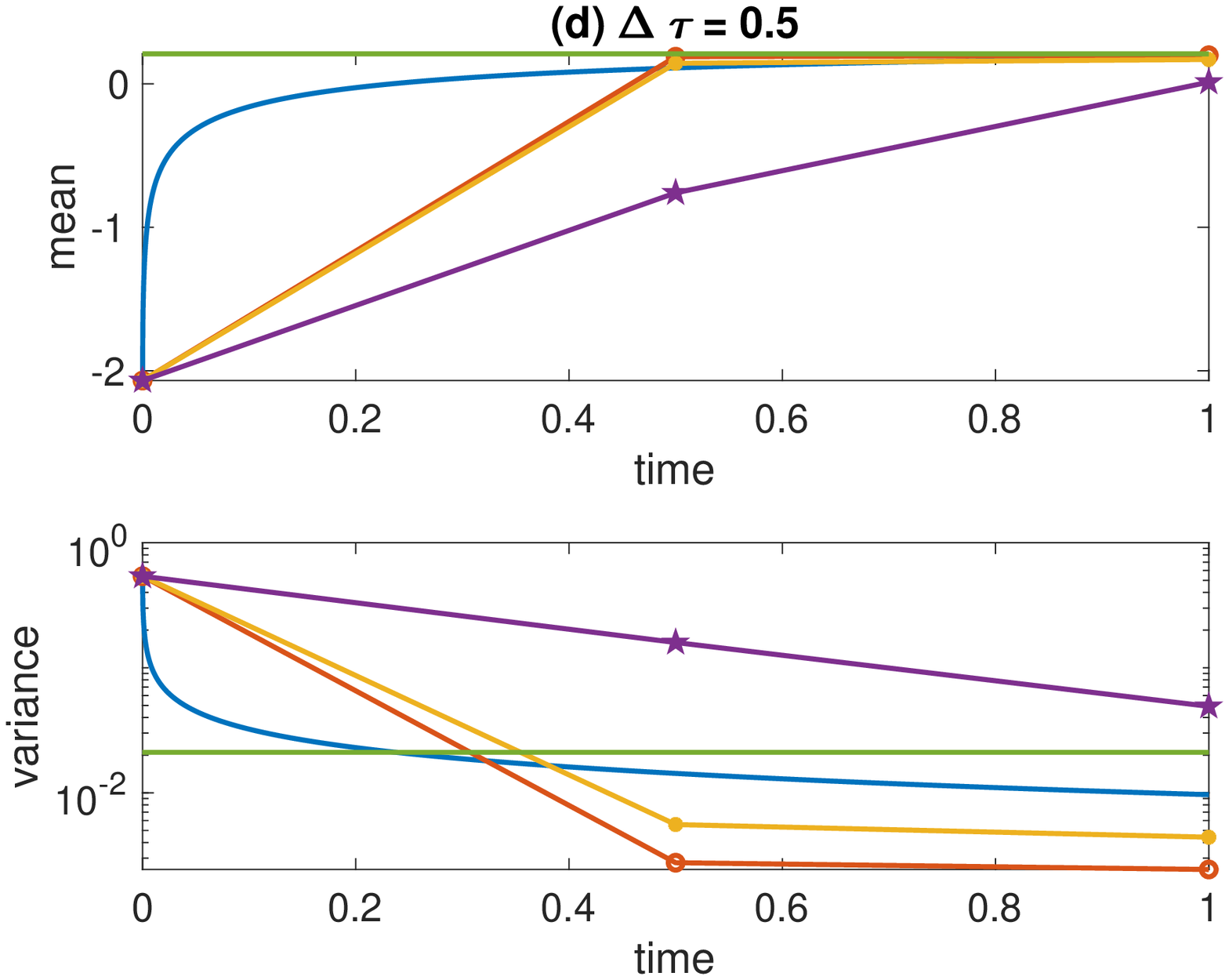} 
\end{center}
\caption{Exact final values and numerical approximations to the mean $m_\tau$ and variance $\sigma_\tau$ estimated by the 
EnKBF using the discrete gradient (DG) with 
$\theta = 1$, the semi-implicit (SI) Euler, and the gradient-free, explicit (IEnKF) method 
with step-sizes (a) $\Delta \tau = 0.01$, (b) $\Delta \tau = 0.1$, (c) $\Delta \tau = 0.2$, (d) $\Delta \tau = 0.5$. The reference time evolution of the mean and 
variance is given by the explicit Euler method with step-size $\Delta \tau = 0.00025$. }
\label{fig:figure2}
\end{figure}

We used $M=100$ particles in our numerical experiments and step-sizes of $\Delta \tau \in \{0.01,0.1,0.2,0.5\}$. 
An explicit Euler discretisation is stable only for the smallest of these step-sizes and is used to provide 
the reference solution. See Figure \ref{fig:figure2} for details.

We find that both the semi-implicit Euler and the discrete gradient method approximate the exact time-evolution of the gradient flow system, as represented
by the explicit Euler method, rather well for step-sizes $\Delta \tau \in \{0.01,0.1\}$. The IEnKF method 
leads to systematically different results in the variance which, however, correspond rather well to the exact posterior value for this particular problem. We also note that the IEnKF method also performs well for the larger step-size of $\Delta \tau = 0.2$. Overall the IEnKF method (\ref{eq:GFE}) 
emerges as the computationally cheapest and 
most accurate method among all the four tested time-stepping procedures.

\subsection{Particle flow Fokker--Planck dynamics}

We repeat the experiments from the previous subsection using particle flow Fokker--Planck dynamics. While its implementation is more
involved, particle flow Fokker--Planck dynamics have the advantage of being more generally applicable than the Kalman-based formulations. 
We also implement the particle flow Fokker--Planck dynamics for a sequential data assimilation procedure for the Lorenz-63 model
\citep{sr:lorenz63}.

\subsubsection{Linear problem}

\begin{figure}
	\begin{center}
		\includegraphics[width=0.45\textwidth,trim = 0 0 0 0,clip]{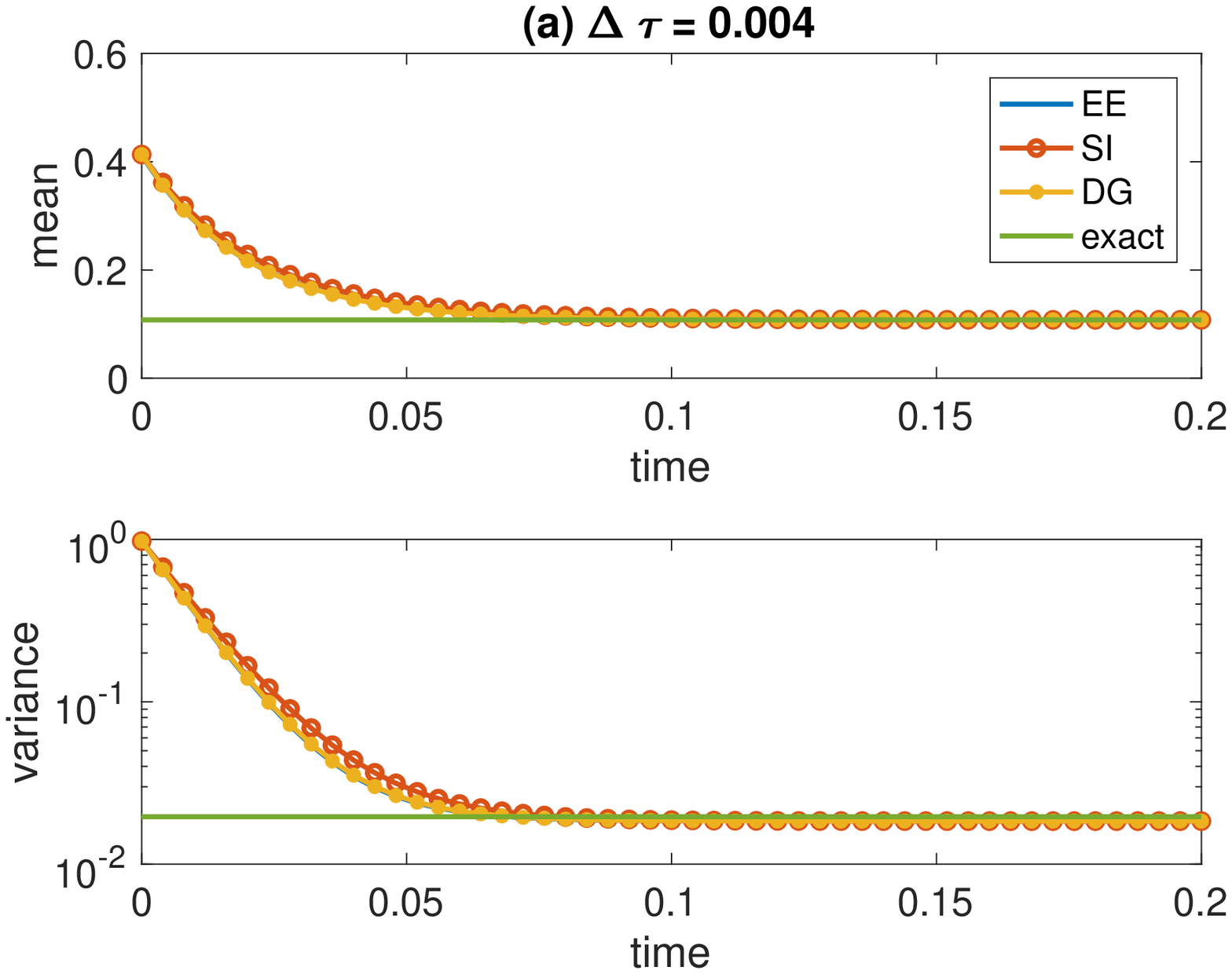} $\qquad$
		\includegraphics[width=0.45\textwidth,trim = 0 0 0 0,clip]{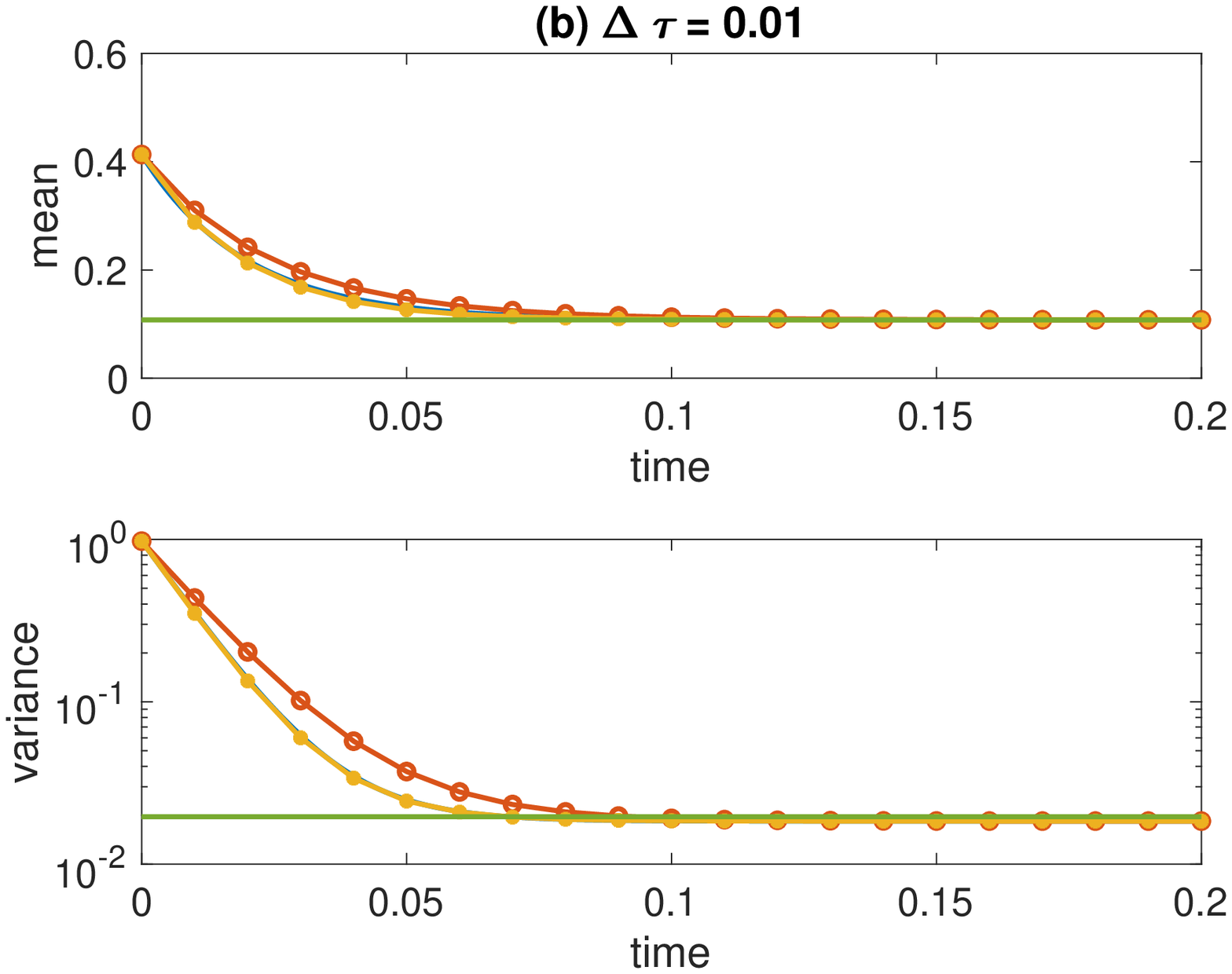} \\ \medskip
		\includegraphics[width=0.45\textwidth,trim = 0 0 0 0,clip]{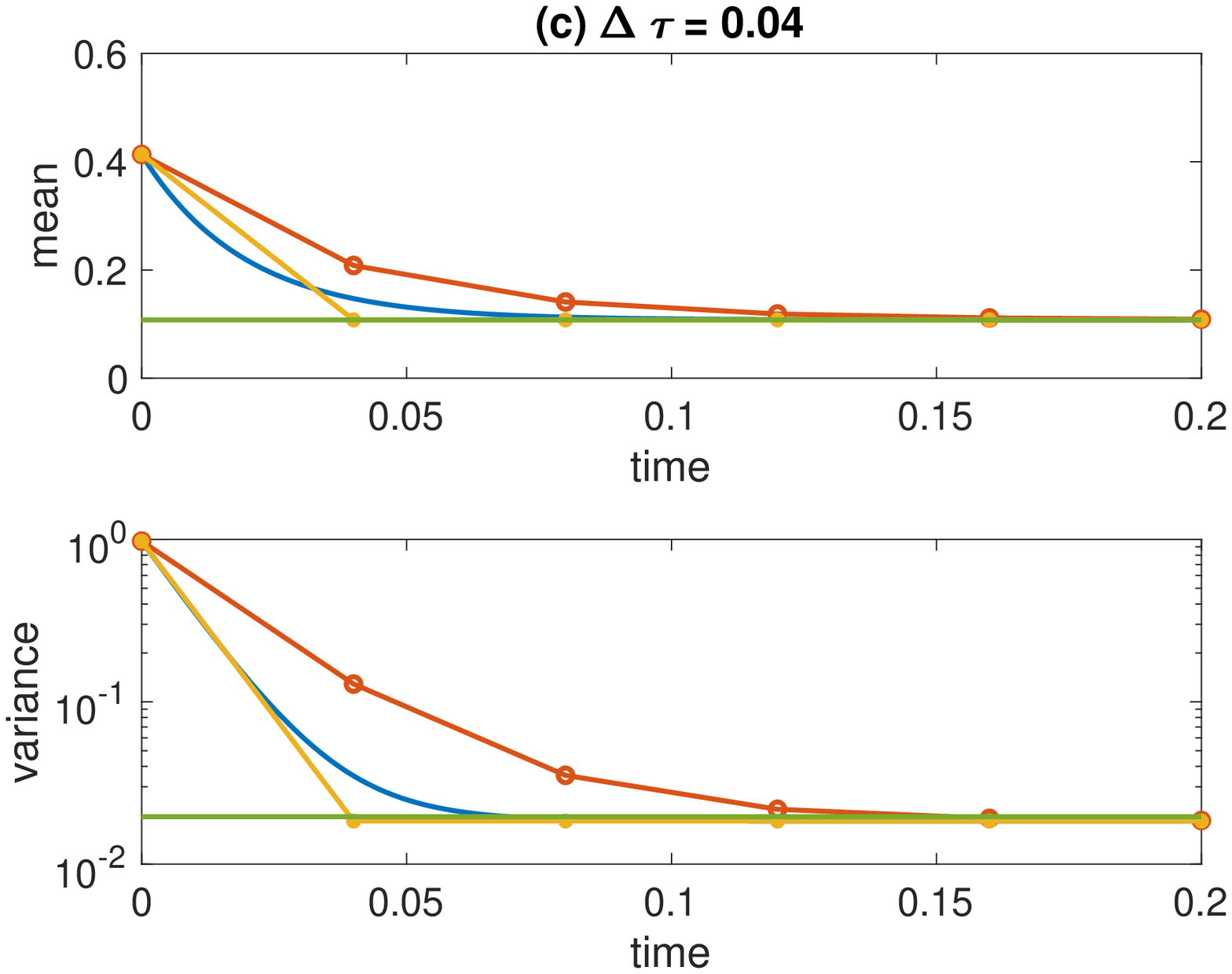} $\qquad$
		\includegraphics[width=0.45\textwidth,trim = 0 0 0 0,clip]{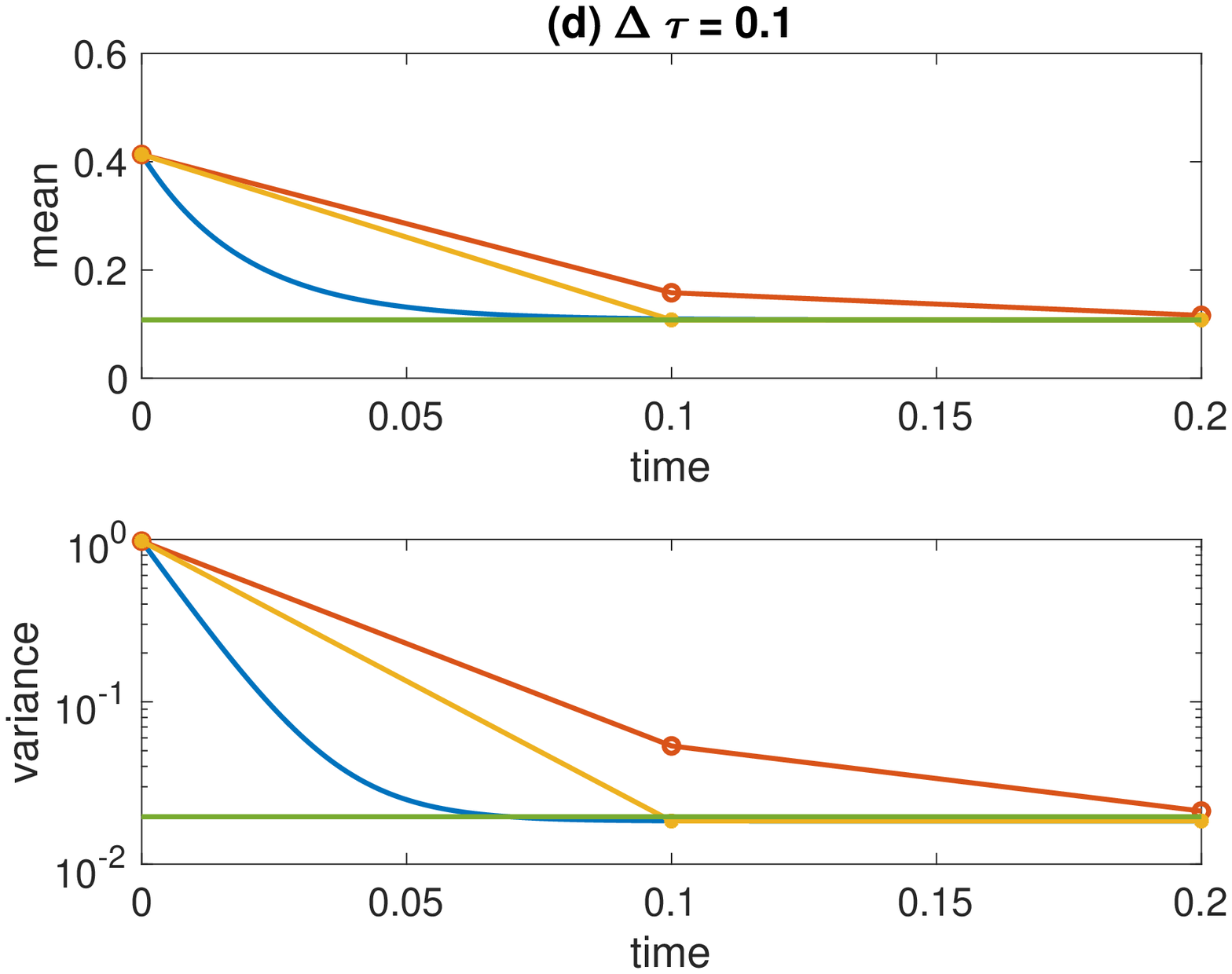} 
	\end{center}
	\caption{Exact final values and numerical approximations to the mean $m_t$ and variance $\sigma_t$ for the linear problem.  Posterior is estimated by the Fokker-Planck dynamics with $M=10$ particles using the discrete gradient (DG) method with 
		$\theta = 1$ and the semi-implicit (SI) Euler method
		with step-sizes (a) $\Delta \tau = 0.004$, (b) $\Delta \tau = 0.01$, (c) $\Delta \tau = 0.04$, (d) $\Delta \tau = 0.1$. The reference time evolution of the mean and 
		variance is given by the Explicit Euler (EE) method with step-size $\Delta \tau = 2 \times 10^{-4}$. }
	\label{fig:figure3}
\end{figure}

We use exactly the same experimental set up as already used for the EnKBF, that is, the measurement error is Gaussian with variance $r = 0.02$
the the prior is Gaussian with mean $m_0 = 1/2$ and variance $\sigma_0 = 1$. We now use $M=10$ particles to approximate the (Gaussian) 
posterior distribution with a Gaussian mixture (\ref{eq:GMFP}). The parameter $\alpha$ in (\ref{eq:shrink}) 
is set to $\alpha = 0.005$. The particle system reaches fairly quickly a stationary configuration which approximates the first two moments of
the posterior distribution rather well. See Figure \ref{fig:figure3}. We also find that the discrete gradient method with $\theta = 1$ relaxes to
the stationary configuration within a single time-step even for large step-sizes $\Delta \tau$.  The semi-implicit Euler method also performs rather well.

\subsubsection{Nonlinear problem}

\begin{figure}
	\begin{center}		
		\includegraphics[width=0.45\textwidth,trim = 0 0 0 0,clip]{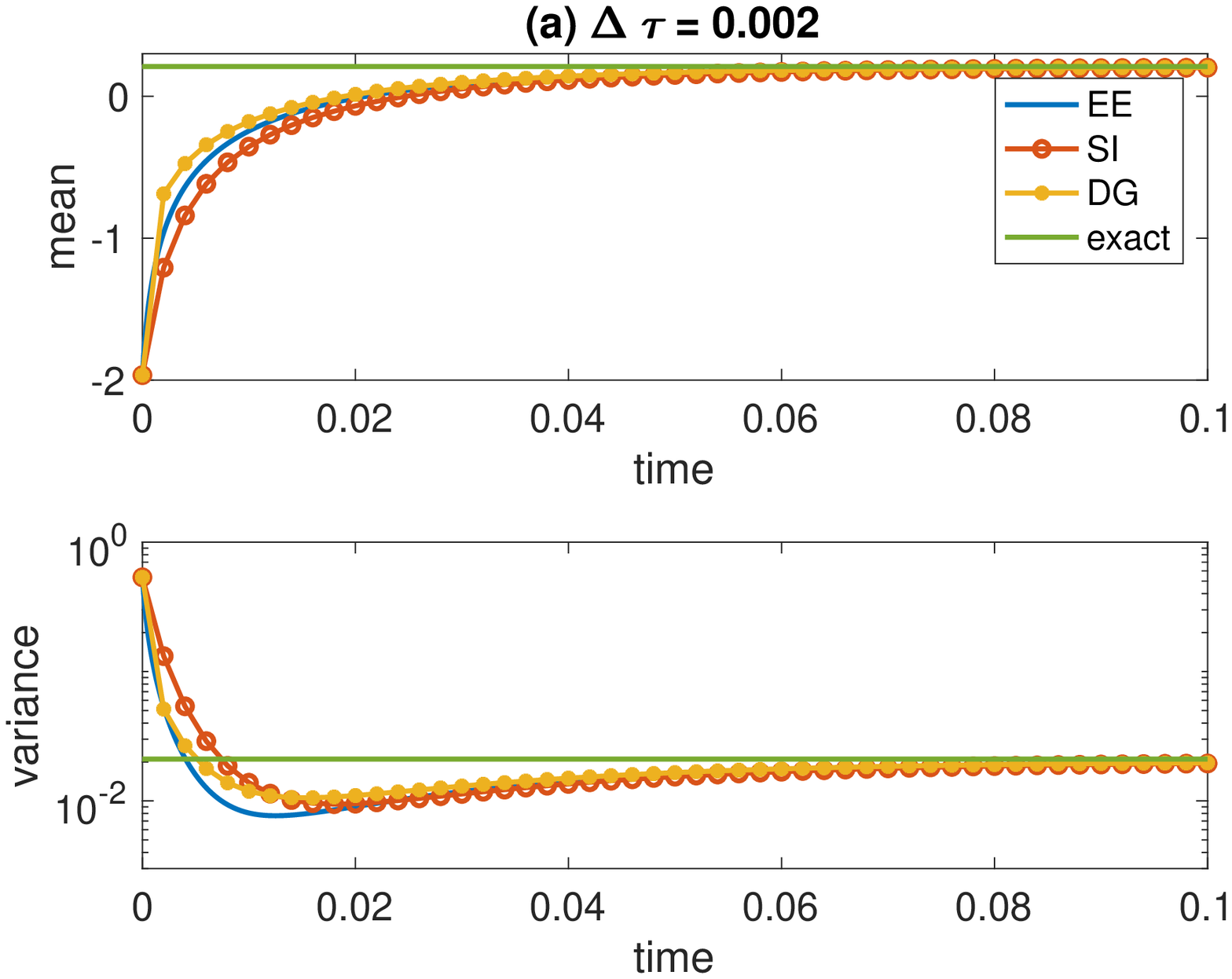} $\qquad$
		\includegraphics[width=0.45\textwidth,trim = 0 0 0 0,clip]{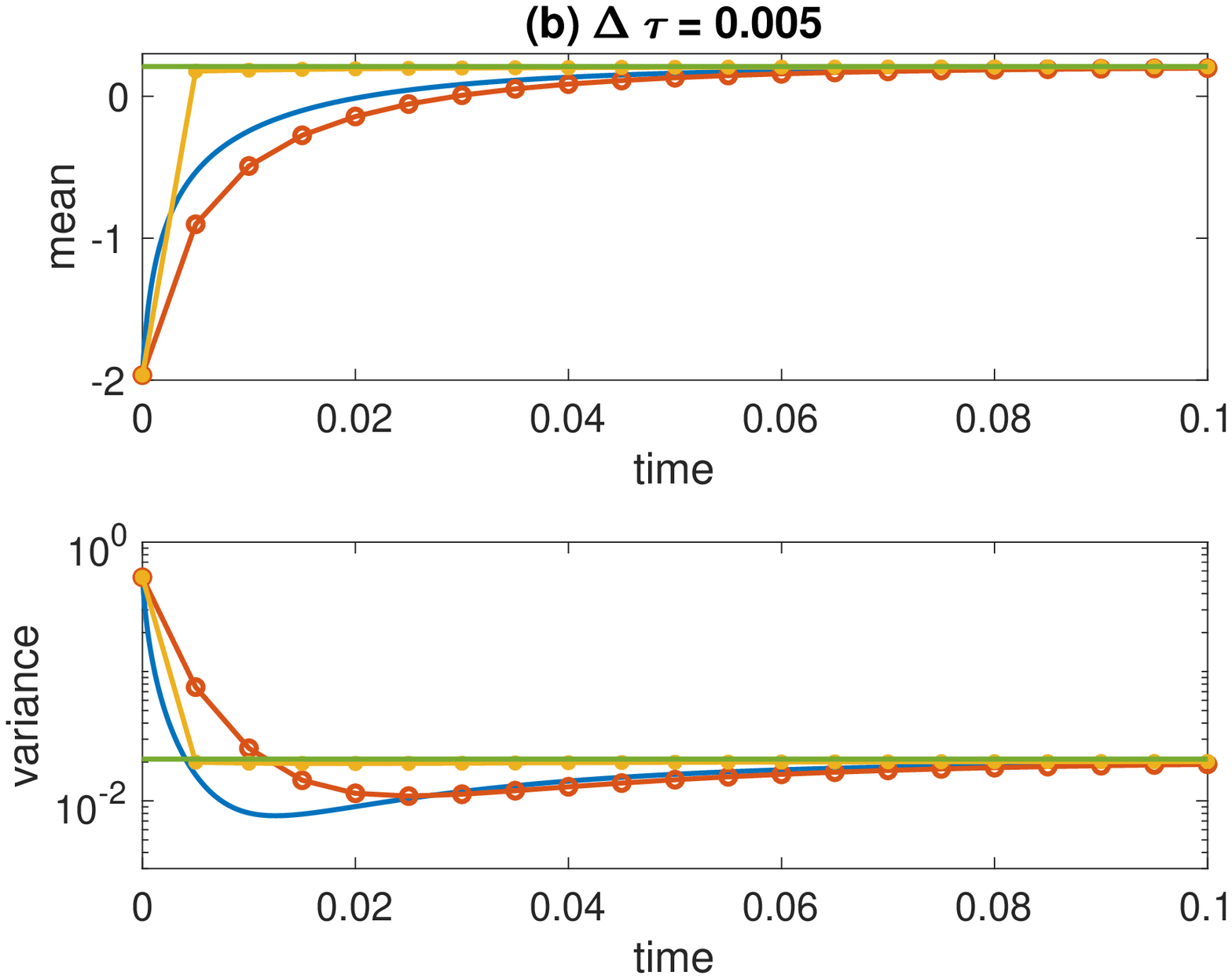} \\ \medskip
		\includegraphics[width=0.45\textwidth,trim = 0 0 0 0,clip]{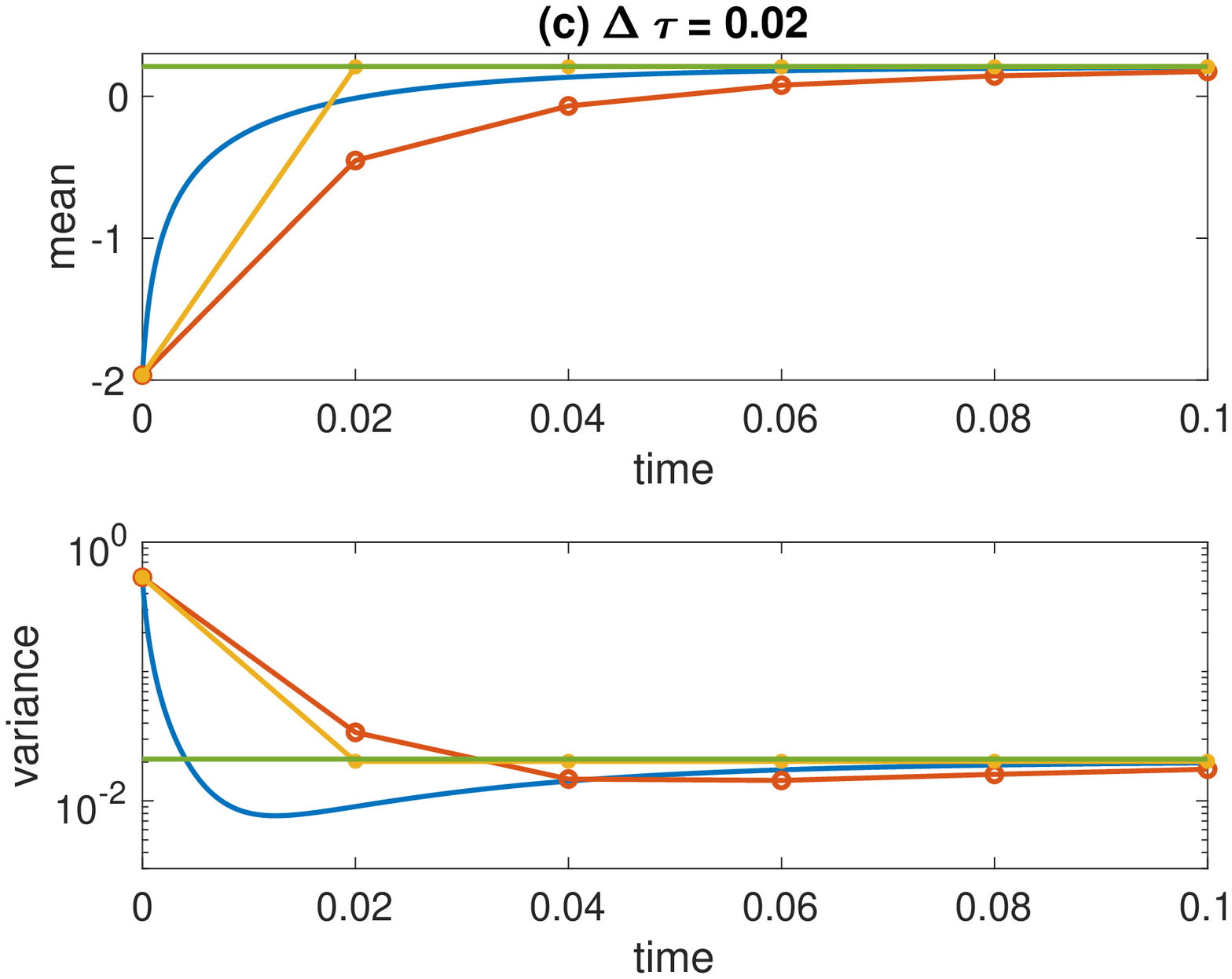} $\qquad$
		\includegraphics[width=0.45\textwidth,trim = 0 0 0 0,clip]{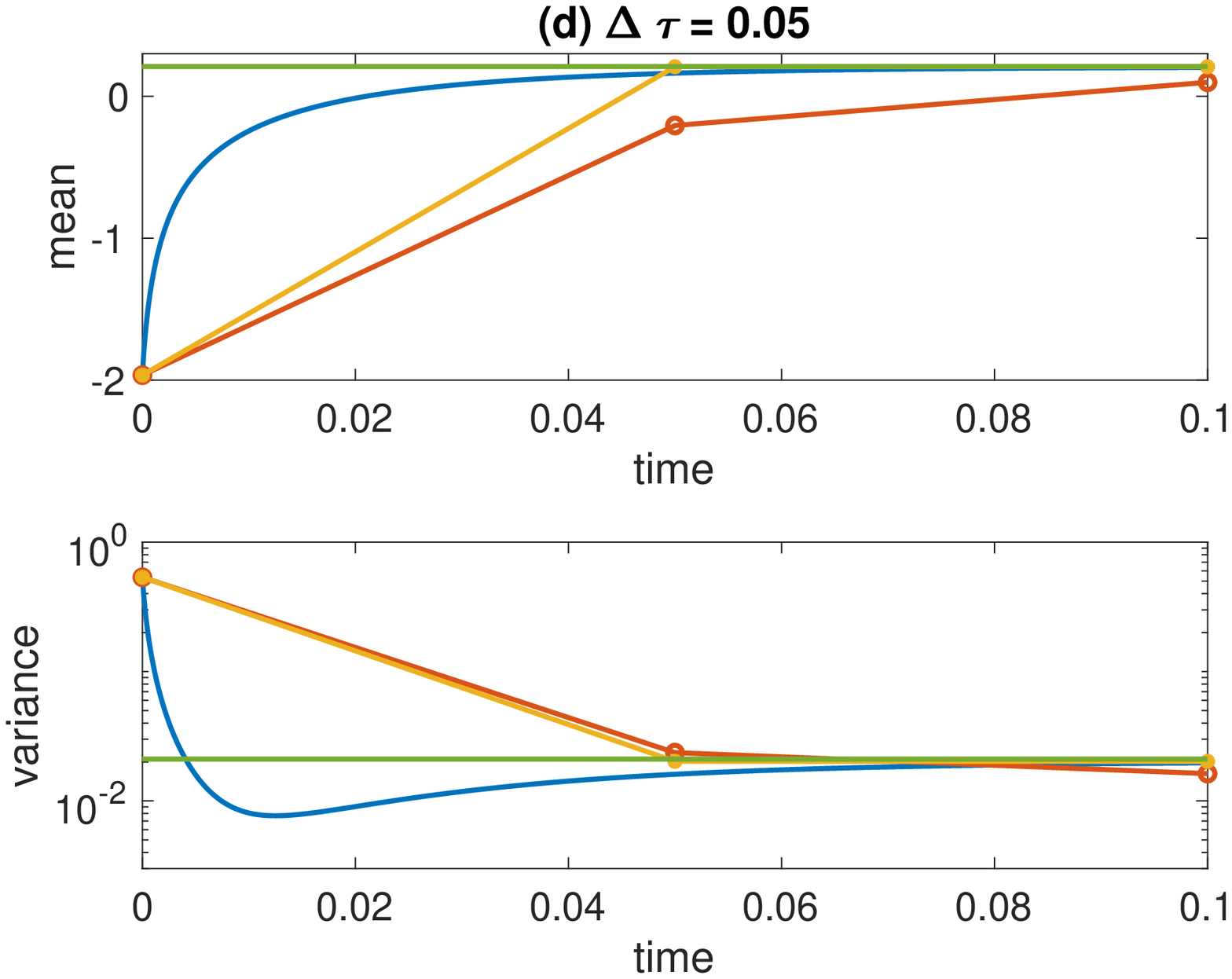} 
	\end{center}
	\caption{Exact final values and numerical approximations to the mean $m_t$ and variance $\sigma_t$ for the non-linear problem.  Posterior is estimated by the Fokker-Planck dynamics using the discrete gradient (DG) method with 
		$\theta = 1$ and the semi-implicit (SI) Euler method
		with step-sizes (a) $\Delta \tau = 0.002$, (b) $\Delta \tau = 0.005$, (c) $\Delta \tau = 0.02$, (d) $\Delta \tau = 0.05$. The reference time evolution of the mean and 
		variance is given by the Explicit Euler (EE) method with step-size $\Delta \tau = 2.5 \times 10^{-6}$. }
	\label{fig:figure4}
\end{figure}

We now come to a more challenging example with a non-Gaussian posterior distribution and 
consider the same nonlinear forward map and experimental setup as described in Section \ref{sec:numnonlinear}.  Again, we use $M = 100$ particles and consider the discrete gradient method with $\theta = 1$ and the semi-implicit Euler method.  The reference solution is also provided by the explicit Euler with small step size 
$\Delta \tau = 2.5 \times 10^{-6}$.  Here we consider step-sizes of $\Delta \tau \in \{0.002, 0.005, 0.02, 0.05 \}$, see Figure \ref{fig:figure4} for results.  The minimisation problem (\ref{eq:MP1}) is solved using the {\tt fminunc} Matlab routine, which is based on a quasi-Newton method. 
The parameter $\alpha$ in (\ref{eq:shrink}) is set to $\alpha = 0.01$.

As for the linear example, we find that the discrete gradient method with $\theta =1$ converges rapidly to the stationary configuration with both the mean and
the variance of the posterior distribution well approximated. The semi-implicit method also behaves well with a only slightly slower convergence rate at
large step-sizes. 

\begin{remark} One could also implement the preconditioned formulation (\ref{eq:DFP1b}) with $\mathcal{A}(\{x_n^j\})$ being the empirical covariance
matrix of the particle system times the ensemble size. This can imply significant computational savings for higher-dimensional problems 
if the kernel matrix $B$ in (\ref{eq:kernel}) is also provided by an
empirical covariance matrix. Furthermore, the explicit computation of gradients can be avoided similar to what is done in the IEnKF implementation (\ref{eq:GFE}). 
More precisely, consider a negative log-likelihood function of the form (\ref{eq:S}), then
\begin{equation}
\nabla_x S(x) \approx \left(P_\tau^{xx}\right)^{-1} P_\tau^{xh} R^{-1}(h(x)-y).
\end{equation}

See also \cite{sr:GHLS19} for related ideas in the context of Brownian dynamics. However, the gradient flow
structure is lost under this approximation.
\end{remark}

\subsubsection{Lorenz-63}

We use the experimental setting of \cite{sr:CRR15}.  Specifically, the Lorenz-63
differential equations \citep{sr:lorenz63} with the standard parameters $\sigma = 10$, $\rho = 28$, and $\beta = 8/3$ are solved numerically with the implicit midpoint rule and step-size $\Delta t = 0.01$. We only observe the first component of the state vector in observation intervals
of $\Delta t_{\rm obs} = 0.12$ with observation error variance $R=8$. A total of $K = 200,000$ assimilation cycles are performed.  The ensemble size varies between $M=15$ and $M=50$ and particle rejuvenation is applied with
$\beta = 0.2$ \citep{sr:CRR15}. The quality of the data assimilation method is assessed using the root mean squared error (RMSE) 
\begin{equation}
\mbox{RMSE} = \frac{1}{K} \sum_{i=1}^K \sqrt{\frac{1}{N_x} \|\overline{x}_{\rm a}(t_k) - x_{\rm ref}(t_k)\|^2}
\end{equation}
with $N_x = 3$ and $t_k = 0.12k$. Here $x_{\rm ref}(t_k) \in \mathbb{R}^3$ denotes the reference solution and $\overline{x}_{\rm a}(t_k)$
the analysis mean at observation time $t_k$. 

At any given time in the sequential data assimilation procedure, a forecast ensemble $\{x_{\rm f}^i\}$ is generated by propagating the posterior sample (also referred to as analysis ensemble) $\{x_{\rm a}^i\}$ from the previous time step through the discretised Lorenz-63 equations.  The transformation of this forecast ensemble to the posterior for the next time step via Gaussian mixture approximations and the particle flow 
Fokker--Planck dynamics is conducted as follows.  
One first computes the empirical mean $\overline{x}_{\rm f}$ and the empirical covariance matrix $P_{\rm f}$ from
the forecast ensemble. Given a parameter value $\alpha \in (0,1]$, a Gaussian mixture approximation to the 
forecast distribution is provided by
\begin{equation}
\pi_{\rm f}(x) = \frac{1}{M}\sum_{i=1}^M {\rm n}(x;\overline{x}_{\rm f}^i,B_{\rm f})
\end{equation}
with covariance matrix $B_{\rm f} = (2\alpha-\alpha^2)P_{\rm f}$ and
\begin{equation} \label{eq:shrink2}
\overline{x}_{\rm f}^i = x_{\rm f}^i - \alpha(x_{\rm f}^i-\overline{x}_{\rm f}),
\end{equation}
as already discussed in Section \ref{sec:FPD}. Given an observation model of the form (\ref{eq:obs}) with linear forward
map $h(x) = Hx$, the analysis (or posterior) distribution is also a Gaussian mixture 
\begin{equation} \label{eq:L63_posterior}
\pi_{\rm a}(x) = \sum_{i=1}^M w^i \,{\rm n}(x;\overline{x}_{\rm a}^i,B_{\rm a})
\end{equation}
with normalised weights
\begin{equation}
w^i \propto \exp \left(-\frac{1}{2} (H\overline{x}_{\rm f}^i-y)^{\rm T}(HB_{\rm f}H^{\rm T} + R)^{-1}(H\overline{x}_{\rm f}^i)\right),
\end{equation}
new centers
\begin{equation}
\overline{x}_{\rm a}^i = \overline{x}_{\rm f}^i - K(H\overline{x}_{\rm f}^i - y), \qquad K = B_{\rm f}H^{\rm T}(HB_{\rm f}H^{\rm T}+R)^{-1},
\end{equation}
and new covariance matrix
\begin{equation}
B_{\rm a} = B_{\rm f} - KHB_{\rm f}
\end{equation}
We next use the particle flow Fokker--Planck dynamics (\ref{eq:DFP1}) with $\pi_\ast = \pi_{\rm a}$, kernel functions
$\psi(x) = {\rm n}(x;0,B_{\rm a})$ and initial particle positions $x_0^i = \overline{x}_{\rm a}^i$ in order to find an equally weighted
Gaussian mixture approximation to (\ref{eq:L63_posterior}). Let us denote the resulting ensemble at a suitably chosen final time $\tau_{\rm end}$
by $\{x_\ast^i\}$. Then the analysis ensemble is finally given by
\begin{equation}
x_{\rm a}^i = x_\ast^i +  B_{\rm a}^{1/2} B_{\rm f}^{-1/2}(x_{\rm f}^i - \overline{x}_{\rm f}^i).
\end{equation}
Note that $\alpha=1$ leads back to the standard ensemble Kalman filter \citep{sr:evensen,sr:reichcotter15}, 
while one formally approaches a standard particle filter for $\alpha \to 0$.
However, it should be noted that the particle flow Fokker--Planck dynamics is not well-defined for $\alpha \to 0$ and fixed ensemble size $M$.

\begin{table}
\begin{center}
\begin{tabular}{|l|c|c|c|c|c|c|c|c|}
\hline
$\alpha\backslash M$ 		&       $15$      & $20$       & $25$      &  $30$   &  $35$    &    $40$      &   $45$        &    $50$    \\
\hline
$0.8$ 	& $2.5145$  &  $2.4726$ &   $2.4471$  &  $2.4406$  &  $2.4221$  &  $2.4173$  &  $2.4162$  &  $2.4112$ \\
$0.85$ & $2.5048$   & $2.4585$ &   $2.4357$  &   $2.4330$  &  $2.4186$  &   $2.4161$  &   $2.4020$  &  $2.4026$ \\
$0.9$  &  $2.4989$  &  $2.4493$ &   $2.4396$  &  $2.4259$  &  $2.4177$  &  $2.4103$ &   $2.4106$ &   $2.4075$ \\
$0.95$ & $2.4998$  &  $2.4605$ &   $2.4399$ &   $2.4327$ &   $2.4242$ &   $2.4132$ &   $2.4162$ &   $2.4131$ \\
$1.0$  & $2.4804$   & $2.4597$ &   $2.4518$  &  $2.4354$  &  $2.4334$  &  $2.4247$  &  $2.4244$  &   $2.4203$ \\
\hline
\end{tabular}
\end{center}
\caption{RMSE for sequential data assimilation applied to the Lorenz-63 model. The parameter $\alpha = 1$ corresponds to
a standard square root ensemble Kalman filter, while $\alpha<1$ results in a Gaussian mixture approximation to the prior and
posterior distributions. Small improvements over the ensemble Kaman filter can be found for ensemble sizes $M\ge 35$.}
\label{table1}
\end{table}

The RMSEs for ensemble sizes $M\in \{15,\ldots,50\}$ and parameters $\alpha \in \{0.8,\ldots,1\}$ in (\ref{eq:shrink2})
can be found in Table \ref{table1}. One finds that $\alpha < 1$ leads to some improvements especially for larger ensemble sizes. Overall
these improvements are, however, not as impressive as those obtained by the hybrid methods considered by \cite{sr:CRR15}. 

%
%
\section{Summary} \label{sec:Summary}
%
%

We have discussed two instances of gradient dynamics arising from Bayesian inference. The gradient flow structure lends itself naturally
to discrete gradient time-stepping method. However, such methods are implicit and potentially costly to implement.
We found that discrete gradient dynamics does not improve the behaviour of iterated ensemble Kalman filters such
as the derivative-free formulation (\ref{eq:GFE}). The situation is different for the particle flow Fokker--Planck dynamics, where the 
discrete gradient method was shown to converge rapidly. This difference in behaviour can be explained by the fact that ensemble Kalman filter
requires an accurate representation of its dynamics while the Fokker--Planck dynamics only requires a rapid transition to an equilibrium solution. 
Applications of particle flow Fokker--Planck dynamics 
to a Lorenz-63 sequential data assimilation problem showed small improvements over a standard ensemble Kalman filter. Further applications
of these methods may arise in the context of machine learning. See, for example, \cite{sr:O17,sr:KS18}.

\medskip

\noindent
{\bf Acknowledgement.} This research has been partially funded by 
Deutsche Forschungsgemeinschaft (DFG) through grant
CRC 1294 \lq Data Assimilation\rq \,\,(projects A02 and B04).

%
%
\bibliographystyle{agsm}
\bibliography{bib_paper}
%
%

\end{document}